\documentclass[12pt]{article}
\usepackage{amsmath}
\usepackage{amsthm}
\usepackage{epic}
\usepackage{epsfig}
\usepackage{graphicx}
\usepackage{floatflt}
\usepackage{doc}
\usepackage{makeidx}
\usepackage{subeqnarray}
\usepackage{stmaryrd}

\catcode`\@=11

\def\ps@headings{\let\@mkboth\markboth
    \def\@oddfoot{}\def\@evenfoot{}
    \def\@evenhead{{\rm \thepage}\hfil { \rm \leftmark}}
    \def\@oddhead{{\rm \rightmark}\hfil {\rm \thepage}}
    \def\chaptermark##1{\markboth{\ifnum \c@secnumdepth >\m@ne
          \@chapapp\ \thechapter. \ \fi ##1}{}}%
    \def\sectionmark##1{\markright{\ifnum \c@secnumdepth >\z@
       \thesection. \ \fi ##1}}}

\ps@headings

\catcode`\@=12

\newtheorem{proposition}{Proposition}[section]
\newtheorem{definition}[proposition]{Definition}
\newtheorem{lemma}[proposition]{Lemma}

\newtheorem{cor}[proposition]{Corollary}

\newcounter{rm}[equation]
\newcommand {\CC}{{\cal C}}

\newcommand {\CG}{{\cal G}}

\newcommand {\CO}{{\cal O}}
\newcommand {\CP}{{\cal P}}
\newcommand {\CQ}{{\cal Q}}
\newcommand {\CR}{{\cal R}}

\newdimen\AAdi%
\newbox\AAbo%
\def\AArm{\fam0 }
\def\AAk#1#2{\setbox\AAbo=\hbox{#2}\AAdi=\wd\AAbo\kern#1\AAdi{}}%
\def\AAr#1#2#3{\setbox\AAbo=\hbox{#2}\AAdi=\ht\AAbo\raise#1\AAdi\hbox{#3}
}%
\def\BBb{{\AArm I\!B}}%
\def\BBn{{\AArm I\!N}}%
\def\BBr{{\AArm I\!R}}%
\def\BBz{{\AArm Z\!\!Z}}%
\newcommand{\pardef}{\stackrel{def}{=}}
\newcommand{\wh}{\widehat}

\newcommand{\disp}{\displaystyle}
\newcommand{\wt}{\widetilde}

\newcommand{\N}{\BBn}

\newcommand{\R}{\BBr}

\newcommand{\Z}{\BBz}

\newcommand{\eps}{\varepsilon}

\newcommand{\inte}[1]{\stackrel{\circ}{#1}}
\newcommand{\ul}[1]{\underline{#1}}

\catcode`\é=\active \defé{\'e}
\catcode`\è=\active \defè{\`e}
\catcode`\ê=\active \defê{\^e}
\catcode`\à=\active \defà{\`a}
\catcode`\ù=\active \defù{\`u}
\catcode`\ô=\active \defô{\^o}
\catcode`\î=\active \defî{\^\i}
\catcode`\ç=\active \defç{\c c}

\newcounter{ex}
\newcounter{qe}
\newcommand{\s}{\sigma}

\newcommand{\us}{{us}}

\newcommand{\llb}{\llbracket}
\newcommand{\rrb}{\rrbracket}
\renewcommand{\S}{\Sigma}
\theoremstyle{definition}
\newtheorem{remark}{Remark}

\begin{document}

\title{Invariant manifolds and equilibrium states for non-uniformly hyperbolic horseshoes}
\author{Renaud Leplaideur\footnote{D\'epartement de math\'ematiques, UMR 6205, Universit\'e de Bretagne Occidentale, 29285 Brest Cedex, France} and Isabel Rios\footnote{Instituto de Matem\'atica, Universidade Federal Fluminense, Rua M\'ario Santos Braga s/n, Niter\'oi, RJ 24.020-140, Brasil} \thanks{This work was partially suported by CNRS-CNPq, UBO, PRONEX-Dynamical Systems, FAPERJ and PROPP-UFF.}
}

\maketitle

\begin{abstract}
In this paper we consider horseshoes containing an orbit of homoclinic tangency accumulated by periodic points. We prove a version of the Invariant Manifolds Theorem, construct finite Markov partitions and use them to prove the existence and uniqueness of equilibrium states associated to H\"older continuous potentials. 
\end{abstract}

\section{Introduction and statement of results}
The goal of this paper is to study some dynamical and ergodic properties of a special class of non-uniformly hyperbolic horseshoes. The non-uniform hyperbolicity, for the systems studied here, comes as a consequence of the presence of a single orbit of homoclinic tangency {\em inside} the horseshoe, that is, accumulated by periodic orbits of it. 

For uniformly hyperbolic systems, results as the existence of stable and unstable manifolds and equilibrium states, as the ones we present here, are obtained in abstract, from a general theory that applies to all systems. The existence of the hyperbolic splitting over the compact invariant set under study, together with the uniform rates of expansion and contraction are strongly used. The conjugacy between the system and some subshift in finitely many symbols allows one to have a complete description of almost all the orbits, in a very  wide sense.

In order to extend this theory beyond the uniformly hyperbolic case, one usually considers two settings: partially hyperbolic systems, and non-uniformly hyperbolic systems. In the first case, the lack of hyperbolicity comes from the degeneracy in the rates of expansion and contraction, and some invariant splitting is assumed to be kept. In that case, under some conditions, the invariant manifolds can be shown to exist (see e.g. \cite{pesin-partial-hyper} for a  survey on the subject). In the second case, the existence of invariant manifolds is shown, for instance, in the so-called Pesin theory, for almost all points, according to some measure (see e.g. \cite{Pes} and \cite{fathi-herman-yoccoz}).  Features as the size of the manifolds, in both cases, depend on what are called hyperbolic returns, and it is not possible to say much about them in abstract.

A significant part of the theory to study non-uniformly hyperbolic dynamical systems is based on models and examples.  In our case we present a model, first introduced in \cite{rios}, and prove an appropriate version of the so-called (un)stable manifold theorem, stating some of the properties of the invariant manifolds. We also prove that the dynamical system is semi-conjugated to the full 3-shift. As a consequence of the regularity of the semi-conjugacy, we obtain existence and uniqueness of an {\it equilibrium state} associated to each given H\"older continuous function $\varphi$.

Here, the setting is in some sense  between the classical topological studies of uniformly hyperbolic dynamical systems and the Pesin's Theory. Since homoclinic tangencies take part in the limit set, we do not have partial hyperbolicity or dominated splitting of the tangent space over it. Moreover, one of the goals of this work is to construct invariant measures. Since the development of Pesin's theory  is, itself, based on an invariant measure, it makes no sense to apply it here.

Though it appears natural from the features of the map that it might be conjugated to the full 3-shift, the construction of the conjugacy follows from estimates involving arbitrarily large iterates of the map and its inverse. On the other hand, since the map is not expansive at the limit set, the existence of equilibrium states does not follow easily from the standard arguments of upper semi-continuity of the metric entropy.

We refer the reader to Hirsch, Pugh and Shub's book \cite{hirsch-pugh-shub}, for the classic theory of invariant manifolds, and to Bowen's book  \cite{bowen}, for the basic
theory of equilibrium states in the uniformly hyperbolic case. 
Recall that an equilibrium state associated to a H\"older continuous potential $\varphi$ is an invariant probability measure that  maximizes the metric pressure associated to $\varphi$, see also the next subsection. 

For non-uniformly expanding maps, it was proved in \cite{krerley1} that there are equilibrium states associated to almost constant H\"older continuous potentials. In \cite{arbeito-matheus-oliveira}, the  same is obtained for small random perturbations of such maps.

Notice that the Sinai-Ruelle-Bowen measures (SRB measures) are equilibrium states associated to (minus log of) the Jacobian of the system in the unstable direction. In that direction, there are abundant examples in the literature where those measures are studied for non-uniformly hyperbolic systems; see \cite{alves-bonatti-viana}, \cite{Viana-Bonatti}, \cite{young-benedicks}, \cite{Hu}, among others. In the specific case of internal (heteroclinic) tangencies such as the ones we study,  a class of codimension two bifurcating maps is provided in \cite{Enrich}, for which there exist SRB measures.

About uniqueness of the equilibrium state, we mention that this result, in the context of non-uniformly hyperbolic systems, is far from obvious. For instance, in
the classical Pomeau-Manneville example (see
\cite{pomeau-manneville}), the lack of hyperbolicity of the system yields to the existence of two ergodic equilibrium states for the potential $\varphi$ studied there.
It could be that, for some special potentials, the presence of the homoclinic
tangency would cause the loss of uniqueness of such special measures. It turns out that, in the special
place of the boundary of the uniformly hyperbolic dynamical systems occupied by the maps that we study here,
that sort of ergodic bifurcation does not occur.

\subsection{Statement of  the results}

\begin{figure}
	\centering
		\includegraphics{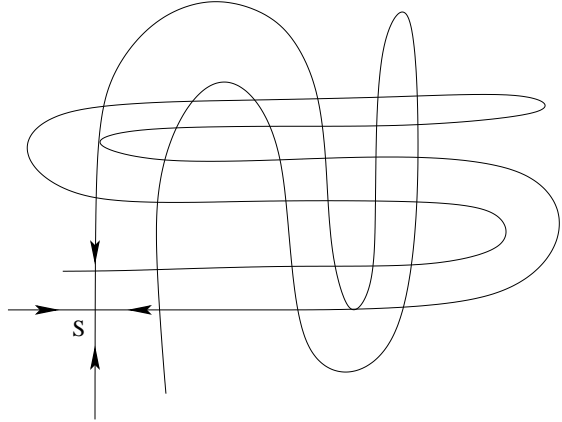}
	\caption{Invariant manifolds}
	\label{bifurcation}
\end{figure}

We study   $\mathcal C ^2$  maps $f$ from the square 
$\CQ=[0,1]\times[0,1]$ into $\R ^2$ with a fixed hyperbolic saddle $S=(0,0)$, whose unstable and stable manifolds have an orbit of homoclinic tangency, as in figure \ref{bifurcation}. 

It was proved in \cite{rios} that, under certain conditions, the set $\Lambda = \cap_{n\in \Z} f^n(\CQ)$ 
admits, for points $x$ outside the tangency orbit, an invariant splitting $T_{x}M= E^s_x \oplus E^u_x$ of the tangent space into stable and unstable directions. Later, in \cite{clrlyap}, the estimate of the rates of contraction and expansion in these 
directions were improved.

Let us assume the notations $\CQ^*=\CQ\setminus \{(0,0)\}$, $\CQ^\#=\CQ\setminus\{(x,y)\in\CQ,\ xy=0\}$, $\Lambda^*=\Lambda\cap \CQ^* \setminus \CO (Q)$, where $\CO (Q)$ is the orbit of homoclinic tangency, and $\Lambda^\#=\Lambda\cap \CQ^\#$. Points in $\Lambda^*$ that are not in $\Lambda^\#$ go to $(0,0)$ for forward or backward iterates.

\medskip \noindent
{\bf Theorem A}
\textit{For every $M$ in $\Lambda^\#$ there exist stable and unstable manifolds  $W^s(M)$ and $W^u(M)$. Moreover the local manifolds $W^s_{l(M)}(M)$ and $W^u_{l(M)}(M)$ are   local graphs from $E^i(M)$ to $E^j(M)$  where $i=u,s$ and $j=s,u$, where $l(M)$
is  continuous  in  $\Lambda^\#$ and tends to 0 when $M$ goes to the critical orbit.
}

\medskip

Let $\varphi:\Lambda \rightarrow \R$ be a H\"older continuous function (we will call such $\varphi$ a {\it potential}). 
We recall that  for any  $f$-invariant probability measure $\mu$, the metric pressure associated to
$\varphi$ is $h_\mu(f)+\int \varphi\,d\mu$, where $h_\mu(f)$ is the
entropy of the measure $\mu$. The quantity $h_\mu(f)+\int
\varphi\,d\mu$ is also referred as the $\varphi$-pressure of the
measure $\mu$. A measure $\mu$ is said to be an \textit{equilibrium
  state} (for $\varphi$) if the $\varphi$-pressure of $\mu$ is maximal
among all $f$-invariant probability measures. In this work we are interested only in probability measures, and any mentioned measure is assumed to be so, from now on.

\medskip \noindent
{\bf Theorem B}
\textit{Given any $\vartheta$-H\"older continuous potential $\varphi$ on
$[0,1]\times[0,1]$, there exists a unique ergodic equilibrium state $\mu_{\varphi}$ for
$f$, associated to the potential $\varphi$. Moreover, $\mu_{\varphi}$ gives positive weight to any open set that intersects $\Lambda$} 

\medskip \noindent

The existence and  uniqueness of the equilibrium state associated to $\varphi$ will follow from a finite-to-one and H\"older continuous semi-conjugacy between the full 3-shift and $\Lambda$.

\subsection{Structure of the paper}

In section \ref{sec-horse} we give the precise definition of the map
$f$ and recall some results about its (non-uniform) hyperbolicity.
In section  \ref{sec-techn-lem} we define a metric adapted to the geometrical features of the system. 
It is also proved there that the map acts in the balls of this metric in a ``Markovian'' way. This fact is
 used later to define a  hyperbolic dynamics $F$, based on $f$ (section \ref{kergodic-charts}), and to prove a version of the invariant manifolds theorem for $F$ and $f$ (section \ref{foliations}). Also in section \ref{foliations}, we  provide some extra properties of the foliations, as the H\"older regularity of the hyperbolic splitting (which is well-known in the uniformly hyperbolic case). In section \ref{sec-thermo} we show that $f$ is not expansive and prove that there exists a H\"older continuous finite-to-one semi-conjugacy from the full 3-shift to $f$. As a consequence, we obtain the existence and uniqueness
of an equilibrium state for any given H\"older continuous potential $\varphi$.

\section{Horseshoes with internal tangencies}\label{sec-horse}
In this section we introduce, in a more precise way, the class of maps that we are going to work with. We define the smooth maps $f$ of the square ${\cal Q}$ into ${\R}^2$,  depending on three parameters, $c$, $\lambda$, and $\sigma$ that admit extensions to the whole plane as smooth diffeomorphisms, having the non-wandering set $\Lambda (f)$ contained in $\cal Q$ (we omit in the notation the dependence on the parameters).

For each allowed choice of the three parameters satisfying some open conditions, the map $f$ is transitive and has a homoclinic tangency associated to the fixed hyperbolic saddle $(0,0)$, which is accumulated by periodic points of the system (we call such a homoclinic point an {\it internal tangency}). These maps where first introduced in \cite{rios}, where some results where obtained for the unfolding of those internal tangencies, and also some properties were found for the map $f$. In \cite{clrlyap} it was shown that the return map to the neighborhood of the internal tangency has nice hyperbolic properties, for instance, the Lyapunov exponents are bounded away from zero.

After introducing the map $f$, also in this section, we recall the construction of the hyperbolic cone fields for points in $\Lambda$ outside the orbit of the tangency (again, we omit the dependence on $f$), and improve some of the estimates on the size of the cones, in order to obtain more accurate bounds for the angles between the stable and unstable directions in each point. 

All of the results here are global, so we need some global control on the region looked at. The computations here use very strongly the definition of the map, we need precise conditions for the first approximation to the problem, but, once it is done, standard methods allow the extension to nearby maps to be done very naturally.

\subsection{The map $f$}

\begin{figure}[htb]
\centerline{\psfig{figure=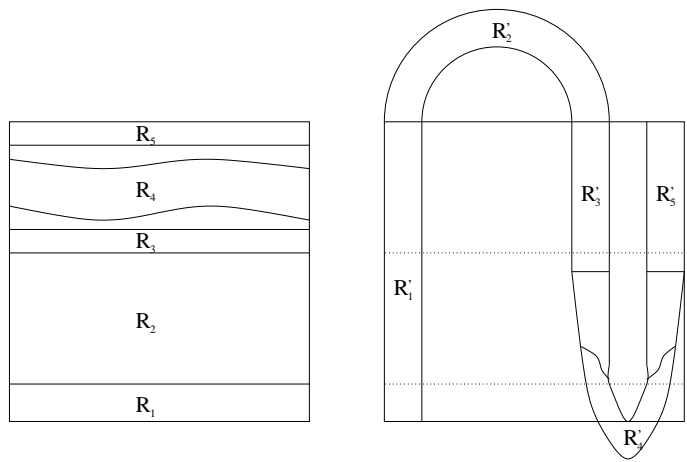,width=9cm}}
\caption{\label{the map}The map $f$}
\end{figure}

Let $\lambda<1/3$, $\sigma>3$. Let $c>0$ be large, some precise conditions on its size are stated along the way. We construct
 a one-to-one differentiable map $f$ from $\cal Q$ into  
${\R} ^2$ satisfying the following conditions (see Figure \ref{the map}):

\begin{itemize}
\item[$a)$] $f(x,y)=(\lambda x,\sigma y)$, if $0\leq y\leq \sigma ^{-1}$ (region $R_1$).
\item[$b)$] $f(x,y)=(\lambda x+(1-\lambda),\sigma y-(\sigma-1))$ if $1-\sigma ^{-1}2/3\leq
 y\leq 1$ (region $R_5$).

\item[$c)$] There exists a horizontal strip, named region $R_3$, contained 
in $[0,1]\times[1/3,1]$, depending on $c$, which is mapped affinely
 in a vertical strip, parallel to the image of the region $R_5$. The
 derivative of $f$ in points of this region is $$Df
(x,y)=\left( \begin{array}{cc}
- \lambda & 0 \\ 0 & -\sigma \end{array} \right).$$
\item[$d)$] Points of $\cal Q$ which are between $R_1$ and $R_3$
 (region $R_2$) are
 mapped outside $\cal Q$.

\item[$e)$] There exists, between $R_3$ and $R_5$,
 a region $R_4$, bounded by two disjoint curves of the form $\{ y=\psi(x):x\in [0,1]\}$, in which the map is not affine, and in this region we have:

\begin{itemize}
\item[$i)$] The top and bottom sides of $R_4$ are mapped into $R_2$,
 outside the image of $R_1$.

\item[$ii)$] $f\left[ \{(0,y):y\in {\R}\}\cap R_4\right]$ is 
contained in the graph of the map $f_0(x)=c(x-q)^2$, with
 $\Vert \frac{\partial f }{ \partial y} (0,y)\Vert \geq \sigma $, where $q \in (2/3,1)$ 

\item[$iii)$] For every $x_0$ in $[0,1]$, $f\left[ \{(x_0,y):y\in {\R}\}\cap R_4\right] $ 
is contained in the graph of the map  $f_{x_0}(x)=c(x-q)^2-\lambda x_0$,
  with  $$\left< \frac{\partial f }{ \partial y} (x,y), 
\frac{\partial f}{ \partial x} (x,y) \right>=0$$  and  
$$\Vert \frac{\partial f }{\partial x}(f^{-1}(q,0))\Vert = \lambda .$$ Notice that we want that the image of $\left[ \{(0,y):y\in {\R}\}\cap R_4\right]$ does not intersect the right side of $\cal Q$.   
\end{itemize}

\item[$f)$] Points between $R_3$ and $R_5$ which are outside $R_4$, are mapped
 inside region $R_2$ with second coordinate greater than $\sigma ^{-1}$. We just 
ask the map to be smooth at this points, and globally one-to-one.
\end{itemize}

In Figure \ref{the map}, $R_i'=f(R_i)$ for $i=1,\ldots ,5$. Notice that $f$ can be extended to ${\R}^2$ in such a way that $(0,0)$ is a hyperbolic fixed point, the left side and the bottom side of $\cal Q$ are contained, respectively, in its unstable and stable manifolds. That implies that $Q=(q,0)$ is a point of homoclinic tangency, whose pre-image we denote by $T=(0,t)$.

These conditions are compatible provided $c$ is big enough, and a precise definition of $f$ in region $R_4$ is given in \cite{rios}. It is also proved there that, for $c$ big enough, a one parameter family of maps that unfolds generically the homoclinic tangency $(q,0)$, crosses the  boundary of the set of uniformly hyperbolic systems at   $f$. That also happens for an open set of nearby families. 

\begin{remark}
\label{composta}
For future uses, we notice that $f_{|R_4}=\Psi\circ \Gamma$, where 
$\Gamma (x,y)=(\lambda x, \sigma (y-t))$ and $ \Psi$ do not depend on $\lambda$ and $\sigma$. 
\end{remark}

\subsection{Hyperbolic cone fields and estimates of angles}

Consider the foliation $\cal F$, of $f(\cal Q)$, whose leaves are images of vertical lines in $\cal Q$ by the map $f$. The leaves of $\cal F$ are vertical in the regions $R_1'$, $R_3'$ and $R_5'$, and parabolic in region $R_4'$. For a point $M$ in this last region, define $E_{\cal P}(M)$ the tangent line at $M$ to the parabola that contains this point.

\begin{remark} \label{rem-isa}
Notice that the biggest value (in modulus) that can be achieved for the slope of $E_{\cal P}(M)$ in the region $R_4'\cap R_1$ is equal to $2\sqrt{c(\lambda + \sigma ^{-1})}$. This slope is achieved at the intersections of the lowest parabolic leaf of $\cal F$ and the line $\{ (x,\sigma ^{-1}):x\in {\R}\}$. The horizontal distance for these two points is $2c^{-1}\sqrt{\lambda + \sigma ^{-1}}$. 
\end{remark}

Due to the linear features of $f$ in most of its domain, the only possible difficulties that can appear to construct a hyperbolic cone field do so close to the tangency points. In order to approach directly this main difficulty, we assume that 
\begin{equation}
\label{tan10}
2\sqrt{c(\lambda + \sigma ^{-1})}<\arctan{\pi /10},
\end{equation} 
increasing $\sigma$ and decreasing $\lambda$, if necessary. Along the way we assume some lower bounds for the value of $c$, and we automatically change $\lambda$ and $\sigma$ in order to keep this property. For reasons to be clarified later, we assume those changes to satisfy $0<b^{-1}<- \frac{\ln{\lambda}}{\ln{\sigma}}<b$, for some positive $b$. We also define 

$$A=R_4'\cap R_1\setminus\{Q\},$$
a region containing a point of the orbit of tangency, that will be looked at in more details throughout this paper.  

We now define a {\em unstable} cone field in $A$, by
assigning for each $M$ in this set a cone ${\cal C}^u(M)$, in the tangent
space to ${\R}^2$ at $M=(x+q,y)$. We recall that the angle $\alpha$
between the horizontal direction and $E_{\cal P}(M)$ satisfies
$\tan{\alpha}=2cx$, and put 

\begin{equation}\label{eqcone}
{\cal C}^u(M)=\left\{(u,v)\in{\R}^2:{\frac{\vert u\vert}{\vert v\vert}}
\leq \frac{\chi_0}{2c|x|}\right\},
\end{equation}
where $\chi_0>1$ is a constant to be precise later. The cone ${\cal C}^u(M)$ is centered at the vertical direction, and, since ${\chi_0}>1$, it contains in its interior the line $E_{\cal P}(M)$. 

\begin{lemma}\label{lem-crois-expo}
Let $n$ be the first positive integer such that $f^n(M)=M'\in A$. Then 
$Df^n_M({\cal C}^u(M))\subset {\cal
C}^u(M')$, and $||Df^n_M(v)||\geq{\sqrt\sigma}^n||v||$, for all $v\in
{\cal C}^u(M)$. If ${\cal
C}^s(M)$ is the
closure of the complement of ${\cal C}^u(M)$,  (the {\em
stable} cone at $M$), we also have
$||Df^{-n}_M(v)||\geq{\sqrt\lambda}^{-n}||v||$, for all $v\in {\cal
C}^s(M)$.

\end{lemma}
\begin{proof} Let us set $M'=(x'+q, y')$. 
In order to prove the lemma, let us do some estimates on the value of $n$. Notice that, before the forward iterates of $M$ are in position to return to $A$, they need to leave the region $[0,1]\times[0,\frac13]$. That gives us 

\begin{equation*}
\sigma^ny\geq\sigma^{n-1}y \geq \frac{1}{3}.
\end{equation*}
Since $y=cx^2- \lambda x_0$, where $0\leq x_0\leq 1$ is the first coordinate of $f^{-1}(M)$, we have 

\begin{equation}
\label{estim1n}
\sigma^ncx^2>\sigma^ny>\frac{1}{3}.
\end{equation}
In order to have $f^{n-1}(M')$ in $A$, we need 

\begin{equation*}
\lambda^{-n+1}x_0' \geq \frac{1}{3},
\end{equation*}
where $0\leq x_0'\leq 1$ is the first coordinate of $f^{-1}(M')$. Since $0\leq y'=c{x'}^2- \lambda x_0'$, we have that
\begin{equation}
\label{estim2n}
\lambda^{-n}c{x'}^2> \lambda^{-n+1}x_0'>\frac{1}{3}
\end{equation}
Those two estimates for the minimum number of iterates give us that

\begin{equation}
\label{estim3n}
n\geq \max \left\{  \frac{\ln{\frac{1}{3c{x'}^2}}}{\ln{\lambda^{-1}}}, \frac{\ln{\frac{1}{3cx^2}}}{\ln{\sigma}} \right\}.
\end{equation}
Now we apply the derivative of $f^n$ at the point $M$ to the vectors of the cone ${\cal C}^u(M)$. Since the map is linear for the $n-1$ first iterates, we have that 

$$Df^{n-1}_M({\cal C}^u(M))=
\left\{(u,v)\in{\R}^2:{\frac{\vert u\vert}{\vert v\vert}}
\leq \left( \frac{\lambda}{\sigma}\right)^{n-1} \frac{\chi_0}{2c|x|}\right\},$$
which is a vertical cone, at the tangent space to ${\R}^2$ at the
point $f^{-1}(M')$. This means, together with the definition of $f$ in
$R_4$, that $Df^{n}_M({\cal C}^u(M))$ is a cone centered at the line
$E_{\cal P}(M')$ such that the oriented angle $\gamma(M')$ between its border lines and $E_{\cal P}(M')$ satisfies 

\begin{equation}
\label{gamma0}
|\tan{\gamma(M')}|< \left( \frac{\lambda}{\sigma}\right)^n \frac{\chi_0}{2c|x|}.
\end{equation}  
Define $\delta (M')=\arctan{\frac{2|x'|c}{\chi_0}}$ (the width of the
stable cone at $M'$, see figure \ref{balle poly}(b)). First we set $\chi_0=4$, such that we have
$\tan{\delta (M')}=1/4 |\tan{\alpha (M')}|$. 

\begin{remark}
For future purposes, we
want to show that $\gamma$ is small enough, not only  to guarantee that
$Df^n_M({\cal C}^u(M))\subset {\cal C}^u(M')$, but also that 

\begin{equation}\label{equ-bonsangles1}
\tan\delta\leq \frac14|\tan\alpha|\leq \frac{3}{4}|\tan{\alpha}|\leq \tan{(|\alpha| -|\gamma|)} \leq \tan{(|\alpha| +|\gamma|)}\leq \frac{6}{5}\tan{|\alpha|},
\end{equation}
and
\begin{equation}\label{equ-bonsangles2}
\tan(|\alpha|-|\gamma|-|\delta|)\geq \frac12 |\tan\alpha|.
\end{equation}
Notice that (\ref{equ-bonsangles1}) and (\ref{equ-bonsangles2}) hold if 

\begin{equation}
\label{equ-bonsangles3}
|\gamma(M')|<\arctan{\frac{6\tan{|\alpha(M')|}}{5}}-|\alpha (M')|.
\end{equation} 
Since we are considering $|\alpha(M)|<\pi/10$ for all $M\in A$ (see (\ref{tan10})), there exists a constant $\chi >0$\index{$\chi$} such that the condition (\ref{equ-bonsangles3}) holds if $\tan{|\gamma}|< \chi \tan{|\alpha|}$. 
\end{remark}

Again, to have this condition, by (\ref{gamma0}) and the choice $\chi_0=4$, it is sufficient to have

\begin{equation*} 
\left( \frac{\lambda}{\sigma}\right)^n \frac{2}{c|x|}< \chi \tan{|\alpha(M')|}=\chi 2c|x'|.
\end{equation*}
that gives us
\begin{equation*} 
\left( \frac{\lambda}{\sigma}\right)^n <\chi c^2|x'||x| .
\end{equation*}
Using (\ref{estim3n}), we find that

\begin{equation}
\label{estimc1} 
\left( \frac{\lambda}{\sigma}\right)^n \leq (3c{x'}^2)^{1-\frac{\ln{\sigma}}{\ln{\lambda}}}
\end{equation}
and

\begin{equation}
\label{estimc2} 
\left( \frac{\lambda}{\sigma}\right)^n \leq (3c{x}^2)^{1-\frac{\ln{\lambda}}{\ln{\sigma}}}
\end{equation}
Considering that $0<b^{-1}<- \frac{\ln{\lambda}}{\ln{\sigma}}<b$, $c$ can be assumed to be $>1$, and  analyzing the two cases $|x|\leq |x'|$ and $|x'|< |x|$, it is enough to have 

\begin{equation}
\label{estimc3} 
(3cx^2)^{1+1/b} < c\chi x^2
\end{equation}
for all $|x|<\lambda <1$. We assume that $c$ is big enough to make valid this relation.

To show that the vectors inside the unstable cone $\CC^u(M)$ grow, by the action of $Df_M^n$, by a factor of at least $\sqrt{\sigma}^n$, just recall that, if $v=(v_1,v_2)\in \CC^u(M)$, then we have
$\Vert v \Vert \leq \sqrt{\frac{4}{c^2x^2}+1}|v_2|$ and $\Vert Df_M^n v \Vert \geq \sigma^n$, giving us

\begin{equation}
\label{estimc4}
\frac{\Vert Df_M^n v \Vert}{\Vert v \Vert} \geq \frac{\sigma^nc|x|}{\sqrt{4+c^2x^2}}=\frac{\sigma^{n/2}(\sigma^{n/2}c|x|)}{\sqrt{4+c^2x^2}}>\frac{\sqrt{c}}{18\sqrt{2}}\sigma^{n/2}>\sigma^{n/2},
\end{equation}
where the last inequalities come from (\ref{estim1n}), the fact
that  $|x|< c^{-1}\sqrt{\lambda + \sigma ^{-1}}<2/c$ and the fact that $c$ can be chosen big. The computations for the
vectors inside the stable cones are analogous; they give another lower
bound (of the same kind) for $c$. 
\end{proof}

Now we extend the unstable cone field
${\CC^u}$ to the whole set $\Lambda$:
for a point $M\in \Lambda$, consider the set
$$I(M)=\left\{n_k \in{\Z}: f^{n_k}(M)\in A \right\}$$
which is the set of ``visits'' of the orbit of $M$ to $A$, where the cone field is already defined.
Based on the fact that 
$$Df^{n_{k+1}-n_k}(M_{n_k})({\CC^u}(M_{n_k}))\subset
{\CC^u}(M_{n_{k+1}}),$$
we choose, for $M_i=f^i(M)$, $n_k<i<n_{k+1}$, cones ${\CC^u}(M_i)$
 such that $Df(M_{i-1})({\CC^u}(M_{i-1}))$ is contained in 
$\inte{{\CC^u}(M_i)}\cup \{ (0,0) \}.$

If $I(M)$ has a first element $n_f$, $M_{n_f}=(x_{n_f},y_{n_f})$, we define, for $i'$ such that $x_{i'}> 1/3$ and $x_i<1/3$ for $i'<i<n_f-1,$
$${\CC^u}(M_{i'})=\left\{(u,v)\in {\R}^2/\frac{\vert {u} \vert }
{\vert {v} \vert}\leq \frac{1}{\sqrt{3}} \right\}.$$ 
Due to the linear features of $f$ at $M_i=(x_i, y_i)$, we have that

$$Df^{n_f-i'-1}{\CC^u}(M_{i'})=\left\{(u,v)\in {\R}^2/\frac{\vert {u} \vert } {\vert {v} \vert}\leq \frac{1}{\sqrt{3}}\frac{\lambda^{n_f-i'-1}}{\sigma^{n_f-i'-1}} \right\},$$
where 
$$\frac{1}{\sqrt{3}}\frac{\lambda^{n_f-i'-1}}{\sigma^{n_f-i'}}\leq 
 \frac{1}{\sqrt{3}}(3c{x_{n_f}}^2)^{1-\frac{\ln{\sigma}}{\ln{\lambda}}}.$$
That, together with (\ref{estimc3}), give us
$$\frac{1}{\sqrt{3}}\frac{\lambda^{n_f-i'-1}}{\sigma^{n_f-i'}}\leq \frac{c\chi x_{n_f}^2}{\sqrt{3}}<\chi 2c|x_{n_f}|,$$
and it is enough to have $Df^{n_f-i'}{\CC^u}(M_{i'})$ included in ${\CC^u}(M_{n_f})$. Now, for $i<i'$, we choose 

$${\CC^u}(M_{i})=\left\{(u,v)\in {\R}^2/\frac{\vert {u} \vert }
{\vert {v} \vert}\leq \frac{1}{\sqrt{3}} \right\},$$ 
and for $i'\leq i <n_f$, we choose cones ${\CC^u}(M_i)$
 such that $Df(M_{i-1})({\CC^u}(M_{i-1}))$ is contained in 
$\inte{{\CC^u}(M_i)}\cup \{ (0,0) \}.$

If $I(M)$ has a last element $n_l$, we take $i'\geq n_l$ the first integer such that $\sigma^{i'-n_l}y_{n_l}>1/3.$ Arguing as before, we set 

$${\CC^u}(M_{i'})=\left\{(u,v)\in {\R}^2/\frac{\vert {u} \vert }
{\vert {v} \vert}\leq \frac{1}{\sqrt{3}} \right\},$$ 
and combine (\ref{estimc3}) with the fact that $c\chi x_{n_l}^2<\sigma^{-1}<1/\sqrt{3}$, give us that $Df^{i'-n_f}{\CC^u}(M_{n_l})$ is included in ${\CC^u}(M_{i'})$. Now, for $n_l<i\leq i'$, we choose cones as before(satisfying the inclusion condition for each iterate), and for $i>i'$, we choose 

$${\CC^u}(M_{i})=\left\{(u,v)\in {\R}^2/\frac{\vert {u} \vert }
{\vert {v} \vert}\leq \frac{1}{\sqrt{3}} \right\},$$ 

If $I(P)$ is empty, we simply set

$${\CC^u}(M_{i})=\left\{(u,v)\in {\R}^2/\frac{\vert {u} \vert }
{\vert {v} \vert}\leq \frac{1}{\sqrt{3}} \right\}.$$ 

That construction provides unstable directions $E^u(P)$ for each point $P$ whose backward iterates are always inside $\CQ$, and stable directions $E^s(P)$ for all points $P$ whose forward iterates are always in $\CQ$. Since the unstable direction $E^u(.)$ is never horizontal, and the stable direction $E^s(.)$ is never vertical, we can fix two unitary vector fields $e^u(.)$ and $e^s(.)$ such that $\left\langle e^u(.),(0,1)\right\rangle >0$ and $\left\langle e^s(.),(1,0)\right\rangle >0$. We have, then, that $Df^{n,-n}_Pe^{u,s}(P)$ are parallel to $e^{u,s}(f^{n,-n}(P))$.

To finish this section, we point out that future changes in the parameter $c$ will keep valid the conditions (\ref{estimc3}) and (\ref{estimc4}), as well as the correspondent ones for the stable case.

\section{Geometric properties of the map $f$}\label{sec-techn-lem}
In this  section we study some geometrical and dynamical features that arise from the definition of the map $f$.  

Let us first state some definitions and notations. We continue to use, for $n\in \BBz$ and $M\in \Lambda$, the notation $M_n$ for the point $f^n(M)$, as in the end of section \ref{sec-horse}. We say that $M$ is in {\it escape phase} if there exists a positive integer $n$ such that $M_{-n} \in A$ and $M_{-i} \in R_1$ for all integer $0\leq i<n$. Analogously,
$M$ is in {\it approach phase} if there exists a positive integer $n$ such that $M_n \in A$ and $M_i \in R_1'$ for all $0< i<n$. If $M=(x,y)$ is in $A$ we set $l(M):=|x-q|$\index{$l$}. If $M$ is in
$\mathcal Q \setminus A$ , we set $l(M):=\sup_{\xi\in A} l(\xi)$. 

Recall that, by definition, the images by the map $f$ of the vertical lines intersected to $A$ are pieces of parabolas that will be called {\it local parabolas} in $A$. If $\CP$ and $\CP'$ are two local parabolas, the closure of the region in $A$ between these two parabolas will be called the {\it parabolic hull} of $\CP$ and $\CP'$.

We denote by $\BBb(M,\eps)$ the ball of center $M$ and with radius $\eps$ for the Euclidean metric $||.||$.
For $M$ in $\Lambda$, let $|v|_M=\max(|v_u|,|v_s|)$, where $v=v_ue^u(M)+v_se^s(M)$. We denote by $B(M,\eps)$ the polygonal ball of center $M$ and radius $\eps$ for this metric $|.|_M$. 

\begin{lemma}
There exists a positive constant $\chi_1$\index{$\chi_1$} such that
for every $M$ in $A$ we have 
\begin{equation}\label{equ-distor-normes}
\chi_1.l(M)|.|_M\leq ||.||\leq 2|.|_M.
\end{equation}
\end{lemma}
\begin{proof}
Let $v$ be any vector in $\R^2$. We set $v=v_ue^u(M)+v_se^s(M)$. Then
we have 
$$||v||\leq |v_u|.||e^u(M)||+|v_s|.||e^s(M)||\leq 2|v|_M.$$
Now, we have $e^u(M)=\cos\theta e^s(M)+\sin\theta e^s_\perp(M)$, where
$\theta$ is some real number and $e^s_\perp$ is the unitary vector
perpendicular to $e^s(M)$ which preserves the orientation of
$\R^2$. For convenience we do the case $|v|_M=|v_s|$; the other case
is similar.
Thus, we have 

\begin{eqnarray}
||v||^2&=& (v_s+v_u\cos\theta)^2+(v_u\sin\theta)^2\nonumber\\
&=& v_s^2+v_u^2+2v_sv_u\cos\theta\nonumber\\
&\geq& v_s^2+v_u^2-2|v_u|.|v_s|\cos\theta\nonumber\\
&\geq & v_s^2\left(
\left(\frac{|v^u|}{|v_s|}\right)^2-2\frac{|v^u|}{|v_s|}.\cos\theta+1\right)\nonumber\\
&\geq& |v|_M^2\sin^2\theta.\label{majo-norm-080704}
\end{eqnarray}
Due to the definition of $A$, we have
$|\sin\theta| \geq \frac{\sqrt 3}{2}|\theta|$. This yields  to the lower bound
in (\ref{equ-distor-normes}) for some constant $\chi_1$.
\end{proof}

For $M$ in $A$, define its escape time as the minimal positive integer $n$ such that $f^n(M)$ does not belong to $R_1$, and the approach time as the minimal positive integer $n$ such that  $f^{-n}(M)\notin R_1'$. 

\begin{remark}\label{lem-2=1+1}
Let $M$ be in $A$, $n_1$ and $n_2$ be its escape and approach times. Then, following the proof of lemma \ref{lem-crois-expo}, we have that $\sigma^{n_1-1}cl^2(M)\geq \frac13$ and $\lambda^{-n_2+1}cl^2(M)\geq \frac13$. An important consequence of this fact is that, for every $0<\rho\leq
1$, the length of the polygonal balls $f^{n_1}(B(M,\rho l(M)))$ and $f^{-n_2}(B(M,\rho l(M)))$, respectively in the
vertical and horizontal direction is at least $\frac\rho4$. 
\end{remark}

\begin{figure}[htb]
\begin{center}
\includegraphics[width=6in]{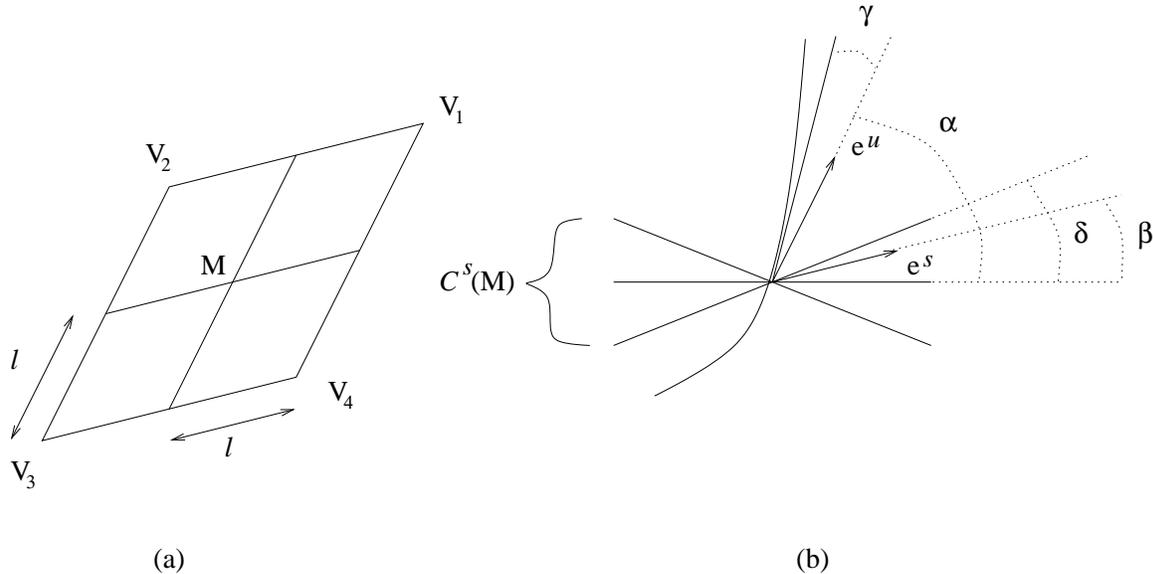}	
	\caption{	\label{balle poly}polygonal ball and angles}
\end{center}
\end{figure}

Figure \ref{balle poly}(a) represents a polygonal ball centered in a point $M\in A$ at the right side of the line $x=q$. Notice that the right and left sides of the polygonal ball are parallel to $e^u(M)$. If the radius $l$ is small, these sides approximate the tangent of the local parabola which contains $M$. Let $S_0(M,l)$ and $S_2(M,l)$ be the bottom and top sides of $B(M,l)$, respectively, and $S_1(M,l)$ be the intersection of $B(M,l)$ with the line parallel to  $E^s(M)$ through $M$. For $S_i(M,l)$ as before, and $0<\rho\leq1$, let $S_i(M,l)(\rho)\subset S_i(M,l)$ denote the segment with radius $\rho.l$ and the same center as $S_i(M,l)$.

\begin{proposition}\label{prop-terrib-para}
There exists a positive real constant $C_0\leq1$\index{$C_0$} such that for every
$M$ in $A$ and for every $0<\rho\leq 1$, if we set $\wt
l (M,\rho)=\rho C_0.l(M)$\index{$\wt l (M)$}, then, every local parabola which crosses the segment $S_1(M, \wt l (M,\rho))(\frac14)$ crosses the two segments $S_0(M, \wt l (M,\rho))$ and $S_2 (M, \wt l (M,\rho))$.
\end{proposition}

\begin{proof}
For simplicity, we assume the notations $\wt l (M,\rho)=\wt l$ and $S_i(M, \wt l (M,\rho))=S_i$. Let $V_1, \ldots ,V_4$ be the vertices of $B(M,\wt l)$ named counter-clockwise from the top right one, as in figure \ref{balle poly}.
Recall that $\alpha$ is the angle between $E_\CP(M)$ and the horizontal, and $\gamma$ is
the angle between the vector $e^u(M)$ and
$E_\CP(M)$ (see figure \ref{polig.eps}(b)). Due to symmetry, we 
assume that $M$ is at the right side of the line $x=q$, without loss of generality.
We denote by $\beta$ the angle between the segments $S_i$ and the
horizontal.
For simplicity we set $q=0$.
In this new system of coordinates, the family of local parabolas  have
equation on the form $Y=cX^2-k$. Let $K$ be such that, for $k=K$, the
associated parabola is the local parabola which contains $M=(x,y)$.
Notice that, though $\alpha$ is positive, it cannot be assumed for $\beta$.
Considering these notations, the coordinates of the points $V_i$ are the following:

$$V_1=\left\{\begin{array}{l}
x_1=x+\wt l(\cos\beta+\cos(\alpha+\gamma))\\
y_1=y+\wt l(\sin\beta+\sin(\alpha+\gamma))\\\end{array}\right.,\ 
V_2=\left\{\begin{array}{l}
x_2=x-\wt l(\cos\beta-\cos(\alpha+\gamma))\\
y_2=y+\wt l(\sin(\alpha+\gamma)-\sin\beta)\\\end{array}\right.$$
$$ 
V_3=\left\{\begin{array}{l}
x_3=x-\wt l(\cos\beta+\cos(\alpha+\gamma))\\
y_3=y-\wt l(\sin\beta+\sin(\alpha+\gamma))\\\end{array}\right.,\
V_4=\left\{\begin{array}{l}
x_4=x+\wt l(\cos\beta-\cos(\alpha+\gamma))\\
y_4=y-\wt l(\sin(\alpha+\gamma)-\sin\beta)\\\end{array}\right.$$ 
For  $0<k\leq 1$, let $V_+(k)$ and $V_-(k)$ be defined by

$$V_+(k)=\left\{\begin{array}{l}
x_+=x+k.\wt l\cos\beta\\
y_+=y+k.\wt l\sin\beta\\\end{array}\right.,\ 
V_-(k)=\left\{\begin{array}{l}
x_-(k)=x-k.\wt l\cos\beta\\
y_-(k)=y-k.\wt l\sin\beta\\\end{array}\right.$$

In this case we have $[V_-(k),V_+(k)]=S_1(k)$. We denote by $\CP^+_k$ and $\CP^-_k$ the two local parabolas containing respectively  $V_+(k)$ and $V_-(k)$. We are looking for sufficient conditions on $\wt l$ and $k$ such that any local parabola which cuts the segment $S_1(k)$ also cuts the segments $S_2=[V_2,V_1]$ and $S_0=[V_3,V_4]$. For this to be satisfied, it is sufficient to have that $\CP^+_k$ and $\CP^-_k$ cut $S_0$ and $S_2$. 

\bigskip\noindent
{\bf - Sufficient conditions  for $\CP^+_k \cap [V_2,V_1]\neq \emptyset$: }

\medskip\noindent
The parabola $\CP^+_k$ crosses $[V_2,V_1]$ if and only if
$$
\left\{\begin{array}{rcl}
y_1&\leq&cx_1^2-2cxk\wt l\cos\beta-ck^2\wt l^2\cos^2\beta+k\wt l\sin\beta-K,\\
y_2&\geq&cx_2^2-2cxk\wt l\cos\beta-ck^2\wt l^2\cos^2\beta+k\wt l\sin\beta-K.\\
       \end{array}
\right.
$$
This yields to  the following system:

\begin{subeqnarray}
\wt l(\sin(\alpha+\gamma)+\sin\beta)&\leq& 2c\wt
lx(\cos\beta+\cos(\alpha+\gamma))+c\wt
l^2(\cos\beta+\cos(\alpha+\gamma))^2\nonumber\\
&&-2cxk\wt l\cos\beta-ck^2\wt l^2\cos^2\beta+k\wt l\sin\beta,
\slabel{equ1.1-terri-para}\\
\wt l(\sin(\alpha+\gamma)-\sin\beta)&\geq& -2c\wt
lx(\cos\beta-\cos(\alpha+\gamma))+c\wt
l^2(\cos\beta-\cos(\alpha+\gamma))^2\nonumber\\
&&-2cxk\wt l\cos\beta-ck^2\wt l^2\cos^2\beta+k\wt l\sin\beta.
\slabel{equ1.2-terri-para}
\end{subeqnarray}
Because $2cx=\tan\alpha$, (\ref{equ1.1-terri-para}) is equivalent to 

\begin{eqnarray}\label{equ-rajout1-050704}\sin(\alpha+\gamma)+(1-k)\sin\beta&\leq &
  \tan\alpha((1-k)\cos\beta+\cos(\alpha+\gamma))\nonumber\\
  && +c\wt l[(\cos(\alpha+\gamma)+\cos\beta)^2-k^2\cos^2\beta].
\end{eqnarray}
Let us assume that $k\leq 1$. Then the second term in the right side
of (\ref{equ-rajout1-050704}) is positive. Moreover
(\ref{equ-bonsangles1}) implies
$\tan\alpha-\tan\beta\geq\frac34\tan\alpha$. Therefore, it is sufficient to
have 

$$\cos(\alpha+\gamma)\frac{\tan(\alpha+\gamma)-\tan\alpha}{\cos\beta}\leq
\frac34(1-k)\tan\alpha,$$ to get (\ref{equ1.1-terri-para}). Again
(\ref{equ-bonsangles1}) together with the fact that
$\cos\beta\geq\frac12$ imply that (\ref{equ1.1-terri-para}) is satisfied if $\disp k\leq\frac{7}{15}$.
From now till the end of this step we assume that $k$ satisfies this
last condition.

Equation (\ref{equ1.2-terri-para}) is equivalent to 

$$c\wt l[\frac{(\cos\beta-\cos(\alpha+\gamma))^2}{\cos\beta}-k^2\cos\beta]-\cos(\alpha+\gamma)\frac{\tan(\alpha+\gamma)-\tan\alpha}{\cos\beta}\leq(1+k)(\tan\alpha-\tan\beta).$$
Now we have $\cos\beta\geq \cos(\alpha+\gamma)$; using again
inequalities (\ref{equ-bonsangles1}), (\ref{equ1.2-terri-para}) holds
if 
\begin{equation}\label{equ-rajout2-050704}
c\wt l\leq \frac34\tan\alpha-1\times2\times\frac14\tan\alpha.
\end{equation}
Hence, (\ref{equ1.2-terri-para}) is satisfied if 
$\disp\wt l\leq \frac14l(M)$.

\bigskip\noindent
{\bf - Sufficient conditions  for $\CP^+_k \cap [V_3,V_4]\neq \emptyset$: }

\medskip\noindent
The parabola $\CP^+_k$ crosses $[V_3,V_4]$ if and only if
$$
\left\{\begin{array}{rcl}
y_3&\geq&cx_3^2-2cxk\wt l\cos\beta-ck^2\wt l^2\cos^2\beta+k\wt l\sin\beta-K,\\
y_4&\leq&cx_4^2-2cxk\wt l\cos\beta-ck^2\wt l^2\cos^2\beta+k\wt l\sin\beta-K.\\
       \end{array}
\right.
$$
This yields to the following system:
\begin{subeqnarray}
-\wt l(\sin(\alpha+\gamma)+\sin\beta)&\geq& -2c\wt
lx(\cos\beta+\cos(\alpha+\gamma))+c\wt
l^2(\cos\beta+\cos(\alpha+\gamma))^2\nonumber\\
&&-2cxk\wt l\cos\beta-ck^2l^2\cos^2\beta+k\wt l\sin\beta,
\slabel{equ2.1-terri-para}\\
-\wt l(\sin(\alpha+\gamma)-\sin\beta)&\leq& 2c\wt
-lx(\cos\beta-\cos(\alpha+\gamma))+c\wt
-l^2(\cos\beta-\cos(\alpha+\gamma))^2\nonumber\\
&&-2cxk\wt l\cos\beta-ck^2\wt l^2\cos^2\beta+k\wt l\sin\beta.
\slabel{equ2.2-terri-para}
\end{subeqnarray}
Using the same kind of computations than just above, equation (\ref{equ2.1-terri-para})  holds if $\disp\wt l\leq \frac{7}{40}l(M)$.\\
Equation (\ref{equ2.2-terri-para}) is equivalent to 
$$c\wt
l[k^2\cos\beta-\frac{(\cos\beta-\cos(\alpha+\gamma))^2}{\cos\beta}]\leq(1-k)(\tan\alpha-\tan\beta)+\cos(\alpha+\gamma)\frac{\tan(\alpha+\gamma)-\tan\alpha}{\cos\beta}.$$
Assuming, moreover, that $\disp k\leq\frac14$, it is sufficient to have $\disp\wt
l\leq \frac12 l(M)$ to get (\ref{equ2.2-terri-para}).\\

\bigskip\noindent
{\bf - Sufficient conditions  for $\CP^-_k \cap [V_2,V_1]\neq \emptyset$: }

\medskip\noindent
Again, $\CP^-_k$ crosses $[V_2,V_1]$ if and only if
$$
\left\{\begin{array}{rcl}
y_1&\leq&cx_1^2+2cxk\wt l\cos\beta-ck^2\wt l^2\cos^2\beta-k\wt l\sin\beta-K,\\
y_2&\geq&cx_2^2+2cxk\wt l\cos\beta-ck^2\wt l^2\cos^2\beta-k\wt l\sin\beta-K.\\
       \end{array}
\right.
$$
This yields to the following system:
\begin{subeqnarray}
\wt l(\sin(\alpha+\gamma)+\sin\beta)&\leq& 2c\wt
lx(\cos\beta+\cos(\alpha+\gamma))+c\wt
l^2(\cos\beta+\cos(\alpha+\gamma))^2\nonumber\\
&&+2cxk\wt l\cos\beta-ck^2\wt l^2\cos^2\beta-k\wt l\sin\beta,
\slabel{equ3.1-terri-para}\\
\wt l(\sin(\alpha+\gamma)-\sin\beta)&\geq& -2c\wt
lx(\cos\beta-\cos(\alpha+\gamma))+c\wt
l^2(\cos\beta-\cos(\alpha+\gamma))^2\nonumber\\
&&+2cxk\wt l\cos\beta-ck^2\wt l^2\cos^2\beta-k\wt l\sin\beta.
\slabel{equ3.2-terri-para}
\end{subeqnarray}
Analogously, equation (\ref{equ3.1-terri-para}) holds if  $\disp\wt l\leq \frac{56}{5}(M)$.\\
In the same way, equation  (\ref{equ3.2-terri-para}) holds if  $\disp\wt l\leq \frac{1}{8}l(M)$.\\

\bigskip\noindent
{\bf - Sufficient conditions  for $\CP^-_k \cap [V_3,V_4]\neq \emptyset$: }

\medskip\noindent
Finally, $\CP^-_k$ crosses $[V_3,V_4]$ if and only if
$$
\left\{\begin{array}{rcl}
y_3&\geq&cx_3^2+2cxk\wt l\cos\beta-ck^2\wt l^2\cos^2\beta-k\wt l\sin\beta-K,\\
y_4&\leq&cx_4^2+2cxk\wt l\cos\beta-ck^2\wt l^2\cos^2\beta-k\wt l\sin\beta-K.\\
       \end{array}
\right.
$$
This yields to the following system:
\begin{subeqnarray}
-\wt l(\sin(\alpha+\gamma)+\sin\beta)&\geq& -2c\wt
lx(\cos\beta+\cos(\alpha+\gamma))+c\wt
l^2(\cos\beta+\cos(\alpha+\gamma))^2\nonumber\\
&&+2cxk\wt l\cos\beta-ck^2\wt l^2\cos^2\beta-k\wt l\sin\beta,
\slabel{equ4.1-terri-para}\\
-\wt l(\sin(\alpha+\gamma)-\sin\beta)&\leq& 2c\wt
lx(\cos\beta-\cos(\alpha+\gamma))+c\wt
l^2(\cos\beta-\cos(\alpha+\gamma))^2\nonumber\\
&&+2cxk\wt l\cos\beta-ck^2\wt l^2\cos^2\beta-k\wt l\sin\beta.
\slabel{equ4.2-terri-para}
\end{subeqnarray}
Equation (\ref{equ4.1-terri-para}) is equivalent to 

$$c\wt
l[\frac{(\cos\beta+\cos(\alpha+\gamma))^2}{\cos\beta}-k^2\cos\beta]\leq(1-k)(\tan\alpha-\tan\beta)-\cos(\alpha+\gamma)\frac{\tan(\alpha+\gamma)-\tan\alpha}{\cos\beta}.$$
Then it is sufficient to have 
$$4c\wt l\leq \left(\frac{9}{16}-\frac25\right)\tan\alpha,$$ to get
the above inequality. Thus (\ref{equ4.1-terri-para}) holds if 
 $\disp\wt l\leq\frac{13}{160}l(M)$.\\
Analogously, equation (\ref{equ4.2-terri-para}) holds if $\disp \wt l\leq 4l(M)$. Now choose $C_0$ such that the bounds above are valid for $\wt l \leq C_0 l(M)$ for $\rho =1$. Then they also hold for $0<\rho<1$, and the proof is complete.
\end{proof}

We say that the parabolas  $\CP^-$ and  $\CP^+$ $u$-cross the ball $B(M, \wt l (M,\rho))$.

For future uses, we need to estimate the size of the balls for which the local parabolas through its points $u$-cross $B(M, \wt l (M,\rho))$.
Keeping the notations in proposition \ref{prop-terrib-para}, we state the following.

\begin{proposition}\label{prop-para-diam}
There is  $0<\eps_0 \leq 1$, uniform in $A$, such that, if $P\in B(M,\eps_0 \wt l (M,\rho))$, then the local parabola through $P$ $u$-crosses the ball $B(M, \wt l (M,\rho))$.
\end{proposition}

\begin{proof}
Define $N_1, \dots ,N_4$ be the vertices of $B(M,\eps_0 \wt l (M,\rho))$, named analogously to $V_1, \dots ,V_4$. Due to proposition \ref{prop-terrib-para}, it is enough to show that the local parabola through $P$ crosses $S_1(M, \wt l (M,\rho))(1/4)$, whose extremes are $V_-$ and $V_+$. For this, we show that the local parabolas through $N_1,\dots ,N_4$ cross $[V_-,V_+]$. The local parabola through $P$ is contained in the minimum parabolic hull containing $N_1, \dots , N_4$, so it must cross $[V_-,V_+]$ as well.

$$V_+=\left\{\begin{array}{l}
x_+=x+\frac{1}{4}\wt l\cos\beta\\
y_+=y+\frac{1}{4}\wt l\sin\beta\\\end{array}\right.,\ 
V_-=\left\{\begin{array}{l}
x_-=x-\frac{1}{4}\wt l\cos\beta\\
y_-=y-\frac{1}{4}\wt l\sin\beta\\\end{array}\right.$$

$$N_1=\left\{\begin{array}{l}
x_1=x+\eps_0\wt l(\cos\beta+\cos(\alpha+\gamma))\\
y_1=y+\eps_0\wt l(\sin\beta+\sin(\alpha+\gamma))\\\end{array}\right.,\ 
N_2=\left\{\begin{array}{l}
x_2=x-\eps_0\wt l(\cos\beta-\cos(\alpha+\gamma))\\
y_2=y+\eps_0\wt l(\sin(\alpha+\gamma)-\sin\beta)\\\end{array}\right.$$

\bigskip\noindent
{\bf - Sufficient conditions  for $\CP_1 \cap [V_-,V_+]\neq \emptyset$: }

\medskip\noindent
The parabola $\CP_1$ crosses $[V_-,V_+]$ if and only if
$$
\left\{\begin{array}{rcl}
y_+&\leq&cx_+^2-2cx(1/4)\wt l\cos\beta-c(1/4)^2\wt l^2\cos^2\beta+(1/4)\wt l\sin\beta-K_1,\\
y_-&\geq&cx_-^2-2cx(1/4)\wt l\cos\beta-c(1/4)^2\wt l^2\cos^2\beta+(1/4)\wt l\sin\beta-K_1.\\
       \end{array}
\right.
$$
where $K_1=cx_1^2-y_1.$
This yields to  the following system:

\begin{subeqnarray}
(1/4-\eps_0)(\sin\beta)&\leq & (1/4-\eps_0)(\tan\alpha\cos\beta)+\nonumber\\
&& (1/16)c\wt l(\cos^2\beta-16\eps_0^2(\cos\beta+\cos(\alpha+\gamma))^2)\nonumber \\
&& +\eps_0\sin(\alpha+\gamma)-\eps_0\tan\alpha\cos(\alpha+\gamma)
\slabel{equ1.1-para-diam} \\
(1/4+\eps_0)(\tan\alpha\cos\beta)&\geq& (1/4+\eps_0)\sin\beta+\eps_0\sin(\alpha+\gamma)
-\eps_0\tan\alpha\cos(\alpha+\gamma)\nonumber \\
&& +(1/16)c\wt l(\cos^2\beta-16\eps_0^2(\cos\beta+\cos(\alpha+\gamma))^2)
\slabel{equ1.2-para-diam}
\end{subeqnarray}
Since $\cos\beta>1/2$, and picking $\eps_0$ small enough, we have

\begin{equation}\label{para-diam-faim}
\cos^2\beta-16\eps_0^2(\cos\beta+\cos(\alpha+\gamma) \geq \cos^2\beta-16\eps_0^2(1+\cos\beta)^2 \geq 0 
\end{equation}
To get (\ref{equ1.1-para-diam}), it is sufficient to have

\begin{equation}\label{para-diam-sono}
\cos(\alpha+\gamma)\eps_0\left( \frac{\tan\alpha-\tan(\alpha+\gamma)}{\cos\beta}\right)
\leq (1/4-\eps_0)(\tan\alpha-\tan\beta),
\end{equation}
which is true, as soon as $\eps_0<15/92$ (see \ref{equ-bonsangles1}).
To get (\ref{equ1.2-para-diam}), it is sufficient to have 

\begin{equation*}
\left( \frac{3}{16}+\frac{7}{20}\eps_0\right)\tan\alpha \geq \frac{1}{16}c\wt l.
\end{equation*}
Again, see expression (\ref{equ-bonsangles1}), and consider that $\tan\alpha=2cx$ and $\wt l =\rho C_0 l(M)=\rho C_0 x$. It is always true, since $0<C_0,\rho\leq 1$.

The computations for the parabolas through the points $N_2, \dots ,N_4$ are analogous, and give other upper bounds to the size of $\eps_0$.
\end{proof}

\begin{figure}[htb]
\begin{center}
\includegraphics[width=6in]{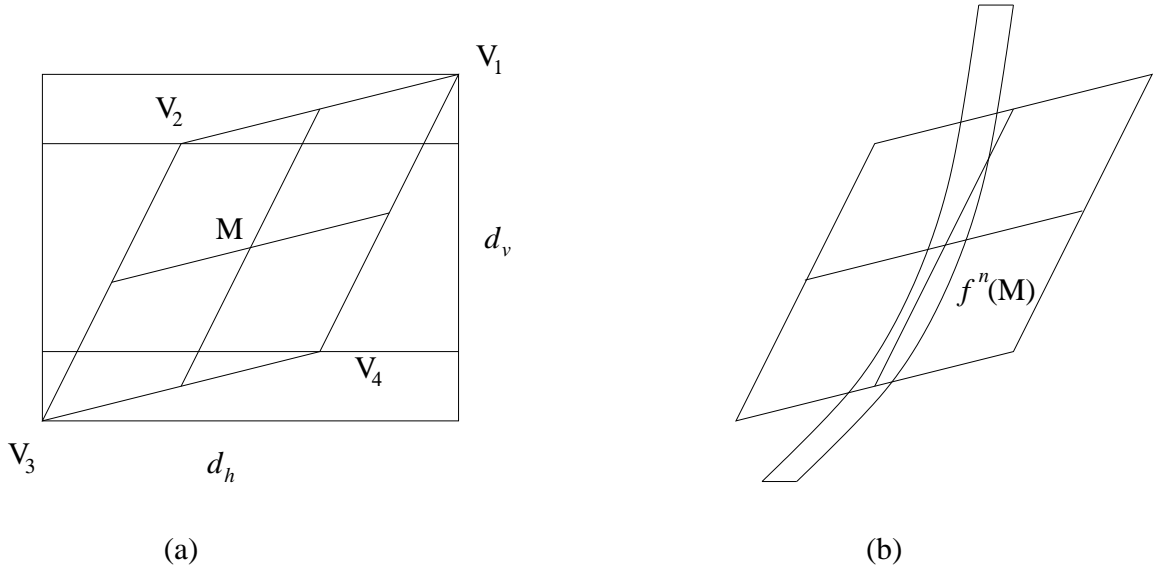}	
	\caption{	\label{distances}polygonal ball and angles}
\end{center}
\end{figure}

Let $M$ be in $A$. For symmetric reasons we
also can assume that $M$ is at the right side of the
homoclinic point $Q$. Consider the ball $B(M,l)$, for some positive $l$, and its vertices $V_i=(x_i, y_i)$ as before. The angle $\alpha$ is positive, but that cannot be assumed for
 $\beta$. Since function $\cos(.)$ decreases close to $0$ and
the function $\sin(.)$
increases close to $0$, then, if $\beta$ is positive, we must have 

$$x_3\leq x_2\leq
x_4\leq x_1\mbox{ and }y_3\leq y_4\leq y_2\leq y_1.$$
Moreover, $|\tan\beta|$ is much smaller than
$\tan\alpha$ and $\tan(\alpha-|\gamma|)$, and thus, if $\beta$ is
 negative we must have 
 
$$x_3\leq x_2\leq
x_4\leq x_1\mbox{ and }y_4\leq y_3\leq y_1\leq y_2.$$
We set $d_v(M, l)=\min(y_1-y_3,y_2-y_4)$ and
$d_h(M,l)=x_1-x_3$. Then $d_h(M,l)$ is the smallest width of a vertical stripe $S_v(M,l)$
containing the
polygonal ball $B(M,l)$, and $d_v(M,l)$ is the largest
height of a horizontal open stripe $S_h(M,l)$ containing $M$ but no edges
of $B(M,l)$. Let $\CR(M,l,l')$ be the rectangle defined by the closure of the
intersection of $S_h(M,l')$ and $S_v(M,l)$.

\begin{definition} If $n$ is the first positive number such that $M_n\in A$, we say that $f^n(\CR(M,l,l'))$ {\it $u$-crosses} the ball $B(M_n,l'')$
if the sides of the connected component of $f^n(\CR(M,l,l'))\cap \CQ$ which contains $M_n$ are parabolas that $u$-cross $B(M_n,l'')$. 
\end{definition}

\begin{proposition}\label{prop-us-cross}
If $\lambda$ and $\sigma^{-1}$ are sufficiently small,
there exists $\eta >0$ such that, for any $M\in A$ and $M_n$ its first
return to $A$, for any $\rho\leq 1$, if $M_n\in B(M',\eta\eps_0\rho C_0l(M'))$ then 
$f^n(\CR(M,\rho C_0l(M),\eps_0\rho C_0l(M))$
 $u$-crosses $B(M',\rho.C_0.l(M'))$ and $B(M',\eps_0\rho.C_0.l(M'))$  (where $\eps_0$ is given in proposition \ref{prop-para-diam}).
\end{proposition}

\begin{proof}

By definition of the map $f$, for $k=0$ to $n-1$, it is linear at $M_k$. 
Consider the notations $\wt l=\rho C_0l(M)$, $\wt l'=\rho C_0l(M')$, $d_v=d_v(M,\wt l)$ and $d_h=d_h(M, \wt l)$. Notice also that $d_v(M,\eps_0\wt l)=\eps_0 d_v$ and $d_h(M,\eps_0 \wt l)=\eps_0 d_h$.  We have $\disp d_v\geq 2\wt
l(\sin(\alpha-|\gamma|)-|\sin\beta|)$. Let us assume, for simplicity
that $\gamma$ and $\beta$ are positive. Applying the mean value theorem to the function $\sin(.)$, and recalling that $\alpha$, $\gamma$ and $\beta$ are small, we have

$$\sin(\alpha-\gamma)-\sin\beta\geq\frac12(\alpha-\gamma-\beta).$$
Moreover, for small positive $t$, $t=\arctan\circ\tan t$. Studying the
derivatives of $t\mapsto \arctan t$, we get
 
$$\alpha-\gamma-\beta\geq \frac12(\tan(\alpha-\gamma)-\tan\beta).$$
Then, (\ref{equ-bonsangles1}) finally yields to
$\disp\sin(\alpha-\gamma)-\sin\beta \geq \frac18\tan\alpha$, and 

\begin{equation}\label{equ1-060704}
d_v\geq \frac12\rho.C_0.cl^2(M).
\end{equation}
By remark \ref{lem-2=1+1}, and the fact that the escape time of $M$ is less than or equal to $n-1$, we find that
$\sigma^{n-1}d_v>\disp \frac16\rho.C_0$. It means that the image by $f^{n-1}$ of the rectangle
$\CR(M,\wt l, \eps_0\wt l)$ is a rectangle with height larger than  $\disp \frac16\eps_0\rho.C_0$
whose center is $M_{n-1}$. The width of this rectangle is
$\lambda^{n-1}d_h$, and we also have $d_h\leq 4\wt l(M)=4\rho.C_0.l(M)$.

Now consider the segment $S_1(M',\rho C_0 l(M'))(1/4)$ and its end points that we call $V_-'$ and $V_+'$. Let $\CP_-'$ and $\CP_+'$ be the local parabolas through these points. Then
$M_{-1}'$ lies between the two pre-images of those curves, which are
vertical lines. We denote by $D_+$ and $D_-$ the horizontal distance
between $M'_{-1}$ and the two vertical lines. 

For each $i=+,-$, the distance $d_i$ equals $\disp
\frac{1}{\lambda}|K(M')-K(M_i')|$, where $K(P)$ is the constance such
that the equation of the local parabola which contains $P$ is
$Y=cX^2-K(P)$. Therefore we get for lower bound

$$D_i\geq \frac{\wt
  l(M')\cos\beta}{4\lambda}|\tan\alpha'-\tan\beta'-\frac{c\wt
  l(M')}{4}\cos\beta'|.$$
Using (\ref{equ-bonsangles1}) and others estimates for the angles, we
get $D_i\geq\disp \frac{5}{32}c.l(M').\wt l(M')$.

Now, notice that, if $M_n\in B(M',\eta\eps_0\rho C_0l(M'))$, we have that the horizontal distance between $M_n$ and $M'$ is less than $2\eta\eps_0\rho C_0l(M')$. This gives us $l(M')\geq l(M_n)- 2\eta\eps_0\rho C_0l(M'),$ and $l(M_n)\geq l(M')- 2\eta\eps_0\rho C_0l(M')$, and we have

\begin{equation}\label{MM'}
\frac{l(M_n)}{2}\leq l(M')\leq 2 l(M_n),
\end{equation}
since we can ask $\eta$ to be smaller than $1/4$. 

Now we must estimate the horizontal distance, $d$, between $M_{n-1}$ and $M_{-1}'$. Assuming the usual notation, let $V_1, \dots ,V_4$ be the vertices of $B(M',\eta\eps_0\rho C_0l(M'))$. Each $V_i$ belongs to a local parabola with equation $y=cx^2-K_i$, and the local parabola through $M'$ has equation $y=cx^2-K'$. Then $d$ is the minimum value of $d_i$, where $d_i=(1/\lambda)|K_i-K'|$, for $i=1,\dots ,4$. We have that

\begin{eqnarray*}
d_i &\leq & \eta\eps_0 \wt l'[\cos\beta'|\tan\alpha'-\tan\beta'|+\cos(\alpha'+\gamma')|\tan\alpha'-\tan(\alpha'+\gamma')|\\
&& +\eta c\wt l'(\cos\beta'+\cos(\alpha'+\gamma'))^2].
\end{eqnarray*}

Now, recall that $\tan\alpha'=2cx'=2cl(M')$. Then (\ref{equ-bonsangles1}) gives us $d_i\leq 10\eps_0\eta\wt l'cl(M')$, if $\eta$ is small enough.
To be sure that $f^n(\CR(M, \wt l,\eps_0 \wt l)$ $u$-crosses $B(M',\eps_0\wt l')$ and  $B(M',\wt l')$, we will put conditions to have 

$$d+\lambda^{n-1}d_h\leq D.$$ 
This fills the conditions for $u$-crossing in the stable direction for both balls, and is achieved if we have

$$4\lambda^{n-1}\wt l +10\eta\eps_0\wt l' cl(M')\leq \frac{5}{32\lambda}\rho.C_0.l^2(M')\eps_0,$$
that means
$$4l(M)\leq \lambda^{-n+1}l^2(M')\eps_0[\frac{5}{32}-10\eta c].$$

By remark \ref{lem-2=1+1}, we have
 that $\lambda^{-n+1}l^2(M_n)\geq 1/3$, and due to (\ref{MM'}), it is
 enough to choose $\eta$ such that $10\eta c<1/32$, and to have

$$4l(M)\leq
 \frac{\eps_0}{3\times 32}\leq \lambda^{-n+1}\eps_0\frac{l^2(M_n)}{4}\frac18\leq \lambda^{-n+1}l^2(M')\left(\frac{5\eps_0}{32}-10\eps_0\eta c \right).$$

For the estimates in the unstable direction, recall that (\ref{equ-distor-normes}) proves that the maximal vertical length
for the ball  $B(M',\rho C_0 l(M'))$ is less than $4l(M')$.
Proposition \ref{prop-terrib-para} proves that the image of $R(M,\wt l,\eps_0 \wt l)$
by $f^n$ has the desired vertical length provided that 
\begin{equation}
  \label{eq:epszero2}
  \sigma^{n-1}\frac{\eps_0}{2}\rho.C_0.cl^2(M)\geq  4\rho.C_0.l(M_n)+4\eta.\eps_0.\rho.C_0.l(M_n).
\end{equation}
As $\eps_0\leq 1$ and we also need some $\eta\leq 1$,the previous inequality is achieved if we have 
$$16l(M')\leq\frac{\eps_0}{6}\leq\frac{\eps_0\sigma^{n-1}cl^2(M)}{2},$$ 
recall remark \ref{lem-2=1+1}. 

Both conditions will hold if $l(\xi)$ is small enough for
every $\xi$ in $A$. Remark \ref{rem-isa} implies that 
$l(\xi)\leq c^{-1}\sqrt{\lambda + \sigma ^{-1}}$, so the result holds 
 for $\lambda$ and $\sigma^{-1}$ small enough.
\end{proof}

\begin{cor}\label{cor-us-cross}
Under the conditions of proposition \ref{prop-us-cross},
 $f^n(R(M, l, l)$  $u$-crosses the ball $B(M_n,\rho.C_0.l(M_n))$.
\end{cor}
 
As a consequence of corollary \ref{cor-us-cross}, we get that the image by $f^n$ of the polygonal ball $B(M,\eps_0 \rho C_0)$ overlaps the ball $B(M_n,\rho C_0)$  in the unstable direction but is ``inside'' in the stable direction, in the sense that it does not cross the lateral sides of $B(M_n,\rho C_0)$.

\section{Kergodic charts}\label{kergodic-charts}
In this section we define a set of local charts, that we call {\it kergodic charts}. These charts are directly inspired on the so-called Lyapunov charts (in Pesin's Theory) and helps us to work with a new map $F$, with good rates of hyperbolicity. Differently from the Lyapunov charts, the kergodic charts preserve some continuity, in the sense that their domains vary continuously with the points. In the next subsections we study the three main ingredients needed to define the charts: the induced map $F$, the polygonal balls, and the control on the distortion of the derivatives.

\subsection{Definition of the map $F$}

The lack of non-uniform hyperbolicity for the system in $\CQ$ is due only to the presence of a single orbit of homoclinic tangency. Away from that orbit, the system is, in some sense, uniformly hyperbolic. 
The map $F$ is defined taking advantage of this fact, as some iterate of the map $f$, depending on the point, in such a way that the hyperbolic estimates hold for the first iterate of $F$. This idea is the heart of the Lyapunov charts. However the Lyapunov charts depend on the point only in a Borel way, not being continuous, in general. The key point in our case is that the first return to the region $A$ is, in fact, uniformly hyperbolic.

Let $F$ be the map defined as follows:
\begin{itemize}
\item[-] if $M$ belongs to $R_3\cup
R_5$, then $F(M)=f(M)$;
\item[-] if $M$ belongs to $R_4\cap (R'_3\cup R'_5)$, then $F(M)=f(M)$;
\item[-] if $M$ belongs to $A$ and $f^n(M)$ belongs to $R_3$ or $R_5$, where $n$ is the escape time for $M$, then $F(M)=f^n(M)$,
\item[-] if $M$ belongs to $A$ and $f^n(M)$ belongs to $R_4$, where $n$ is the escape time for $M$, then $F(M)=f^{n+1}(M)$,
\item[-] if $M$ belongs to $A\cap\left([0,1]\times\{0\}\right)$ then we set $F(M)=(0,0)$. We also set $F((0,0))=(0,0)$
\end{itemize}
Notice that $F$ is not defined for points in $R_1\cap R'_1$, except for $(0,0)$. In fact, the other points in that region belong to orbits of points in the rest of the square, and the results of this section apply to them, by iteration of $f$. We also have that $F$ is not injective only at the pre-image of $(0,0)$, being one-to-one in the rest of the domain.  For convenience we denote by $\CQ_F$, the domain of the map $F$, and use the notations
  $\CQ^*_F=\CQ_F\cap \CQ^*$,  $\CQ^\#_F=\CQ_F\cap \CQ^\#$,  $\Lambda_F=\Lambda\cap \CQ_F$,   $\Lambda^*_F=\Lambda_F\cap \Lambda^*$, and  $\Lambda^\#_F=\Lambda_F\cap \CQ^\#$.

\subsection{Polygonal $us$-balls}

Now we start to define the $us$-balls. The definitions takes into account the number of iterates in the future and in the past that a point needs before it visits the region $A$.

Let $\rho\in(0,1)$. For $M$ in $A$, we denote by $B^\us(\rho)$ the
polygonal ball $B(M,\rho.C_3.l(M))$, where $C_3$ is a constant to be estimated later.
If $M\in\Lambda\setminus A$ is between two visits of its orbit to $A$, then let $k$ and $p$ be the smallest positive integer satisfying
$M_{-k}\in A$ and $M_p\in A$. We set 
$$B^\us(M,\rho)= \left(B(M,\rho.C_3t_u(M))\cap
E^u(M)\right)\times\left(B(M,\rho.C_3t_s(M))\cap E^s(M)\right),$$
where 
$$
t_u(M)\pardef\min(\frac13,l(M_{-k})||df^k_{M_{-k}}.e^u(M_{-k})||),$$
 and $$
t_s(M)\pardef\min(\frac13,l(M_p)||df^{-p}_{M_{p}}.e^s(M_{p})||).$$ 
If $M$ is before the first visit of its orbit to $A$, that is, if $M_{-k}\notin A$ for every
$k\geq 0$ and  there exists a positive integer $p$ such that
$M_p\in A$, then we set 
$$B^\us(M,\rho)=\left(B(M,\frac13\rho.C_3)\cap
E^u(M)\right)\times\left(B(M,\rho.C_3t_s(M))\cap E^s(M)\right),$$
where $t_s(M)$ is defined as above.

As the next case, let $M$ be such that $M_p\notin A$
for every $p\geq 0$ and there exists  a positive $k$ such that $M_{-k}\in A$ ($M$ is after the last visit of its orbit to $A$). Then
we set
$$B^\us(M,\rho)=\left(B(M,\rho.C_3t_u(M))\cap
E^u(M)\right)\times\left(B(M,\frac13\rho.C_3)\cap E^s(M)\right),$$
where $t_u(M)$ is defined as above. 

Finally, if the orbit of $M$ does not visit $A$, that is, if
$\forall k\in \Z$, $M_k\notin A$, we set 
$$B^\us(M,\rho)=B(M,\frac13\rho.C_3).$$
These balls are called the $us$-balls. For convenience, $B^i(M,\rho)$ will denote the one-dimensional ball $B^\us(M,\rho)\cap E^i(M)$ ($i=u,s$).

\subsection{Dynamic for the $us$-balls and  definition of the kergodic charts}

One of the key points to prove  the existence of the invariant manifolds, which is the goal of the next section, is to control the distortion of the map $\wh F$ of the kergodic charts, to be defined in this section. Before that, we point out some of the properties of the map $F$.

\begin{itemize}
\item The map $F$ is locally linear, except in the region $R_4$. Thus, by definition of the $us$-balls, for every $M$ in $R_3\cup R_5$, the image, by the map $F$, of the $us$-ball $B^\us(M,\rho)$ is a polygonal ball which $u$-crosses the polygonal ball $B^\us(F(M), \rho)$;  here $u$-crossing means that $F(B^\us(M,\rho))$ is of the kind $\left(B(F(M),\rho.C_3t'_u(M))\cap
E^u(M)\right)\times\left(B(M,\rho.C_3t'_s(M))\cap E^s(M)\right)$, where $t'_u(F(M))\geq t_u(F(M))$ and $t'_s(F(M))\leq t_s(F(M))$.

\item For the same reason, when $M$ belongs to $A$ and $F(M)$ belongs to $R_3$ or $R_5$,
the image by the map $F$ of the $us$-ball $B^\us(M,\rho)$ is also a
polygonal ball which also $u$-crosses the polygonal ball $B^\us(F(M),
\rho)$. In both cases, due to the linearity of $F$, we have $F=DF$. 

\item Let $M$ be in $A$  such that $F(M)\in A$.
Corollary \ref{cor-us-cross} proves that for every $C_3\leq C_0$, for every $\rho$ in $(0,1)$, the image by $F$ of the polygonal $us$-ball $B^\us(M,\rho)$ $u$-crosses the polygonal $us$-ball $B^\us(F(M),\rho)$. 
\end{itemize}

The map $f$ is $C^2$ in the compact set $R_4$; thus there exists some constant $C_4$\index{$C_4$} such that, for every $\xi$ in $R_4$, $ ||D^2f(\xi)||\leq C_4$ (for the Euclidean norm). Hence, for $M \in A$ with $F(M)\in A$, the relation (\ref{equ-distor-normes}) and the previous discussion prove that there exists some positive constants $C_5=C_5(C_4,\chi_1)$ such that, for every $\xi$  and $\xi'$ in the connected component of $B^\us(M,\rho )\cap F^{-1}(B^\us(F(M),\rho))$ which contains $M$,  

\begin{subeqnarray}\label{equ-majo-dF}
\left|DF_\xi-DF_{\xi'}\right|_{F(M)}&\leq& C_5.l^{-1}(F(M))|F(\xi)-F(\xi')|_{F(M)},\slabel{equ-majo-dF1}\\
\left|DF^{-1}_{F(\xi)}-DF^{-1}_{F(\xi')}\right|_M&\leq& C_5.l^{-1}(M)|\xi-\xi'|_{M}.\slabel{equ-majo-dF2}
\end{subeqnarray}
With some adjustments, $C_5$ can be chosen to be uniform in $\Lambda_F^{\#}$, extending to the other cases.  We can now use these properties to define the kergodic charts. We  define them in the same way that the Lyapunov charts are usually defined  (e.g. in \cite{ledrap1}):

\begin{definition} We call kergodic charts, the family of embedings $\Phi_M:B_M(0,\rho)\subset \R^2\rightarrow \R^2$, for $M\in{\Lambda_{F}^\#}$, satisfying:

\begin{description}
\item[(i)] $\Phi_M$ is affine, $\Phi_M(0)=M$ and $D\Phi_M(0)$ respectively maps $\R^u\pardef\R\times\{0\}$ onto $E^u(M)$ and $\R^s\pardef\{0\}\times\R$ onto $E^s(M)$.
If $v$ belongs to $\R^u\cup \R^s$, then $\left|\Phi_M(v)\right|_M=l(M)||v||$.
\item[(ii)] The set $B_M(0,\rho)$ is defined by $\Phi_M^{-1}(B^\us(M,\rho))$; it is provided with the adapted norm $|\,.\,|$: for $v=v_1.(1,0)+v_2.(0,1)$, $|v|=max(|v_1|,|v_2|)$.
\end{description}
We call kergodic maps, the family of maps $\wh F_M\pardef \Phi_{F(M)}^{-1}\circ F\circ \Phi_M$. The family of maps $\wh F_M^{-1}\pardef \Phi_{F^{-1}(M)}^{-1}\circ F^{-1}\circ \Phi_M$ will be called the inverse kergodic maps.
\end{definition}

We also set $B_M^i(0,\rho)\pardef \R^i\cap B_M(0,\rho)$ ($i=u,s$).
 The kergodic maps and charts satisfy the following properties:

\begin{enumerate}
\item For any $\xi$ and $\xi'$ in $B_M(0,\rho)$.  
$$\frac12||\Phi_M(\xi)-\Phi_M(\xi')||\leq |\xi-\xi'|\leq \frac{1}{\chi_1}l^{-1}(M).||\Phi_M(\xi)-\Phi_M(\xi')||.$$
\item For every $M$ in $\Lambda$:
\begin{itemize}
\item[2.1] For every $v$ in $\R^u$, $\disp {\sqrt \s}^{k}|v|\leq |D\wh F_M(v)|$, where $F(M)=f^k(M)$,
\item[2.2] For every $v$ in $\R^s$, $\disp |D\wh F_M(v)|\leq {\sqrt\lambda}^{k}|v|$, where $F(M)=f^k(M)$.
\end{itemize}
Away from $A$,  when $F\equiv f$, this
follows from the uniform hyperbolicity. In $A$ this follows from lemma
\ref{lem-crois-expo}. 
\item For every $\xi$ and $\xi'$ in the connected component of $B_M(0,\rho)\cap \wh F^{-1}_{F(M)}(B_{F(M)}(0,\rho))$ that contains $0$, we have 
\begin{itemize}
\item[3.1]  $\disp \left| \wh F_M(\xi)-\wh F_M(\xi')-D\wh F_M(0).(\xi-\xi')\right|\leq   C_3.C_5.|D\wh F_M(0)(\xi-\xi')|$,
\item[3.2] $\disp \left| \xi-\xi'-D\wh F^{-1}_{F(M)}(0).(\wh F_M(\xi)-\wh F_M(\xi'))\right|\leq C_3.C_5|D\wh F_{F(M)}^{-1}.(\wh F_M(\xi)-\wh F_M(\xi'))|$.
\end{itemize}
\end{enumerate}
These two inequalities follow from (\ref{equ-majo-dF1}) and (\ref{equ-majo-dF2}).
Indeed, the first one comes from the following estimates:

\begin{eqnarray*}
\left| \wh F_M(\xi)-\wh F_M(\xi')-D\wh F_M(0).(\xi-\xi')\right| & = & \left|\int_0^1(D\wh F_M(t(\xi-\xi'))-D\wh F_M(0))(\xi-\xi')dt\right| \\ & \leq & \frac{C_3C_5l(F(M))}{l(F(M))}|D\wh F_M(0)(\xi-\xi')|.
\end{eqnarray*}
Notice that if $F(M)=f^n(M)$, then $f$ is linear at least in $n-1$ iterates. Hence, the possible distortion is only due to the last iterate, and is uniformly bounded. The second inequality is analogous.

\begin{remark}
Should it be convenient, we consider the balls $B_M(0,\rho)$ as subsets of $\R^2$, keeping the index $M$ to recall that the radius of the balls depend on the point $M$.
\end{remark}

\section{Invariant manifolds and some of their properties}\label{foliations}
In this section we extend the theory of invariant manifolds and local product structure to the non-uniformly hyperbolic set $\Lambda^*$. The invariant manifolds are constructed first for the map $F$, and then extended to $f$ by iteration. Since the hyperbolicity of the splitting degenerates as one approaches the tangency, the radius of the neighborhood in which the product structure holds is not uniform in $\CQ$. In the last subsection we give some extra properties for these manifolds.

\subsection{The invariant manifolds}
\begin{proposition}\label{prop-vect-loc-inte}
There
exists a positive real number $\rho_1\leq 1$, such that, if we set $C_3:=\rho_1.C_0$ then, for every $M$
in ${\Lambda_F^\#}$,  and for any
$0<\rho\leq 1$, 
there exists an uniquely defined pair of  curves, $W^u_{\rho}(M)$ and $W^s_{\rho}(M)$, satisfying

\begin{itemize}
\item[-] $W^{u}_{\rho}(M)$ is tangent to $E^{u}(P)$ for each $P\in W^u_{\rho}(M)\cap_{n \in \N} F^n(\CQ_F)$;
\item[-] $W^{s}_{\rho}(M)$ is tangent to $E^{s}(P)$ for each $P\in W^s_{\rho}(M)\cap_{n \in \N} F^{-n}(\CQ_F)$;
\item[-] $W^{u,s}_{\rho}(M)$ is the graph of a function $g^{u,s}_M$, from  $B^{u,s}(M,\rho)$ to $B^{s,u}(M,\rho)$ with Lipschitz constant smaller than $1/3$; 
\item[-] $F(W^u_{\rho}(M))\supset W^u_{\rho}(F(M))$ and $F^{-1}(W^s_{\rho}(M))\supset W^s_{\rho}(F^{-1}(M))$
\end{itemize}
\end{proposition}

\begin{proof}
The proof follows from the graph transform applied to the kergodic maps. The key point in this classical proof (see e.g. \cite{Hirsch-Pugh}) is to control the Lipschitz-closeness between $\wh F_M$ and the linear map $d\wh F_M(0)$. In the following steps, we show that the main conditions to apply the graph transform are satisfied by  $\wh F_M$, if $\rho$ is small enough.

We do now the estimates in detail for the case where $M\in A$ and $F(M)\in A$. In fact, in the other cases, $F$ is linear or uniformly hyperbolic, and the estimates are analogous in the uniformly hyperbolic case. Let $Lip_1(M)$ denote the space of  $1$-Lipschitz functions $s$, from $B^u_M(0,\rho)$ to $B^s_M(0,\rho)$, satisfying $s(0)=0$. Consider in $Lip_1(M)$ the norm 
$$\disp ||s||= \sup_{x\neq 0}\frac{|s(x)|}{|x|}.$$ 
Denote by $\pi_1$ (resp. $\pi_2$) the projection onto $\R^u$ (resp. $\R^s$) in the direction of $\R^s$ (resp. $E^u$).

\medskip
{\bf Step 1 -} The function $\pi_1\circ \wh F_M\circ (Id,s)$, where $Id$ is the identity at $\R^u$, is a Lipeomorphism 
from $B^u_M(0,\rho)$ onto $B^u_M(0,\rho)$, satisfying 
\begin{equation}
\label{lipeo}
Lip(\pi_1\circ \wh F_M\circ (Id,s))^{-1}\leq \frac{1}{(1-\rho\chi)|D\wh F_M(0)_{|\R^u}|},
\end{equation}
for a constant $\chi$, independent of $M$. Indeed, we apply ($3.1$) to 
\begin{eqnarray*}
\xi & = & (x,s(x)) \\
\xi'& = & (y,s(y))
\end{eqnarray*}
and remember that the vector $(x-y,s(x)-s(y))$ is almost horizontal, due to the fact that $s\in Lip_1(M)$. That gives  $|D\wh F_M(0)(\xi-\xi')|=|D\wh F_M(0)_{|\R^u}|.|x-y|$. Applying the argument  on the Lipschitz constant for the inverse map (see \cite{Hirsch-Pugh}, point $1.5$, page $137$), the result follows, for  $\rho$ small enough (remember that $C_3=\rho. C_0$).

\medskip
{\bf Step 2 -} We now set 
$$\Gamma_F(s)=\pi_2\circ \wh F_M\circ (Id,s)\circ (\pi_1\circ \wh F_M\circ (Id,s))^{-1},$$ 
the classical graph transform operator, and we claim that $\Gamma(Lip_1(M)\subset Lip_1(F(M))$. Indeed, we have $\Gamma(s)(0)=0$, since $\wh F_M(0)=0$, and, for $x$ and $y$ in $B^u_{F(M)}(0,\rho)$, 
\begin{eqnarray*}
|\Gamma_F(s)(x)-\Gamma_F(s)(y)| & = & |\pi_2\circ \wh F_M\circ (Id,s)\circ (\pi_1\circ \wh F_M\circ (Id,s))^{-1}(x)\\
 && -  \pi_2\circ \wh F_M\circ (Id,s)\circ (\pi_1\circ \wh F_M\circ (Id,s))^{-1}(y)|.
\end{eqnarray*}
Now, expression (\ref{equ-majo-dF1}) gives that $D\wh F_M(\xi)\circ D\wh F_{F(M)}^{-1}(0)$ is $\chi\rho$-close to the identity in the canonical basis (in the kergodic charts), where $\chi$ is an universal constant (some other conditions will be given for $\chi$ later). Therefore, for any $(\alpha,\beta)$ in $\R^2$, we have 
$$\pi_1(D\wh F_M(\xi) (\alpha (1,0) +\beta (0,1)))=\alpha\chi\rho D\wh F_M(0)_{|\R^u}.(1,0)+\beta(1+\chi\rho)D\wh F_M(0)_{|\R^s}.(0,1).$$
Hence, we get
\begin{equation}
\label{gamma}
|\Gamma_F(s)(x)-\Gamma_F(s)(y)|
\leq \frac{C_3C_5\rho}{(1-\chi\rho)}\frac{|D\wh F_M(0)_{|\R^u}|}{|D\wh F_M(0)_{|\R^u|}|}|x-y|.
\end{equation}
For $\rho$ sufficiently small, $\Gamma_F(s)$ is in $Lip_1(F(M))$. 

\medskip
{\bf Step 3 -} If $x\in B^u_M(0,\rho)$ and $y\in B^s_M(0,\rho)$, then  

$$|\pi_2\circ\wh F_M(x,y)-\Gamma_F\circ\pi_1\circ \wh F_M(x,y)|\leq \chi|D\wh F_M(0)_{|\R^s}|.|s(x)-y|$$
for some constant $\chi$. In fact, we have
\begin{eqnarray*}
|\pi_2\circ\wh F_M(x,y)-\Gamma_F\circ\pi_1\circ \wh F_M(x,y)| & \leq &  
|\pi_2\circ\wh F_M(x,s(x))-\Gamma_F\circ\pi_1\circ \wh F_M(x,y)|\\
 && +  |\pi_2\circ\wh F_M(x,s(x))-\pi_2\circ\wh F_M(x,y)|.
\end{eqnarray*}
Now, expression (\ref{gamma}) gives 
\begin{equation*}
|\pi_2\circ\wh F_M(x,s(x))-\Gamma_F\circ\pi_1\circ \wh F_M(x,y)|\leq
|\pi_1\circ\wh F_M\circ (Id,s)(x)-\pi_1\circ\wh F_M(x,y)|.
\end{equation*}
Now, 
\begin{equation*}
\wh F_M\circ (Id,s)(x)-\wh F_M(x,y)=\int^1_0D\wh F_M(x(1,0)+t(s(x)-y)(0,1)).|s(x)-y|(0,1)dt.
\end{equation*}
Using again (\ref{equ-majo-dF1}), we find 
\begin{equation*}
|\pi_2\circ\wh F_M(x,s(x))-\Gamma_F\circ\pi_1\circ \wh F_M(x,y)|\leq \chi\rho|D\wh F_M(0)_{|\R^s}|.|s(x)-y|. 
\end{equation*}
The same computations produce 
\begin{equation}
\label{eq3}
|\pi_2\circ\wh F_M(x,s(x))-\pi_2\circ\wh F_M(x,y)|\leq \chi|D\wh F_M(0)_{|\R^s}|.|s(x)-y|,
\end{equation}
and, finally, we get
\begin{equation*}
|\pi_2\circ\wh F_M(x,s(x))-\Gamma_F\circ\pi_1\circ \wh F_M(x,y)|\leq \chi|D\wh F_M(0)_{|\R^s}|.|s(x)-y|,
\end{equation*}
for some positive $\chi$.

\medskip
{\bf Step 4 -} It remains to show that $\Gamma_F$ is a contraction from $Lip_1(M)$ to $Lip_1(F(M))$. For $x$ in $B_{F(M)}(0,\rho)$ and for $s_1$ and $s_2$ in $Lip_1(M)$, we have
\begin{eqnarray*}
|\Gamma_F(s_1)(x)-\Gamma_F(s_2)(x)| & = & |\pi_2\circ \wh F_M\circ(Id,s_1)\circ(\pi_1\circ\wh F_M\circ(Id,s_1))^{-1}(x)\\
&& - \Gamma_F(s_2)[\pi_1\circ \wh F_M\circ(Id,s_1)\circ(\pi_1\circ\wh F_M\circ(Id,s_1))^{-1}(x)]|
\end{eqnarray*}
Therefore (\ref{eq3}) gives
\begin{eqnarray*}
|\Gamma_F(s_1)(x)-\Gamma_F(s_2)(x)| & \leq & \chi|D\wh F_M(0)_{|\R^s}|.|s_1(\pi_1\circ\wh F_M\circ(Id,s_2))^{-1}(x)\\ 
&& -  s_1(\pi_1\circ\wh F_M\circ(Id,s_1))^{-1}(x)|.
\end{eqnarray*}
Using the norm for $s_1-s_2$, we get 
\begin{equation*}
|\Gamma_F(s_1)(x)-\Gamma_F(s_2)(x)|\leq \chi|D\wh F_M(0)_{|\R^s}|.||s_1-s_2||\frac{|\pi_1\circ\wh F_M\circ(Id,s_2))^{-1}(x)|}{|x|},
\end{equation*}
and (\ref{eq3}) yields to 
\begin{equation*}
|\Gamma_F(s_1)(x)-\Gamma_F(s_2)(x)|\leq \frac{\chi|D\wh F_M(0)_{|\R^s}|}{(1-\chi\rho)|D\wh F_M(0)_{|\R^u}|}||s_1-s_2||.|x|,
\end{equation*}
and, finally, 
\begin{eqnarray*}
|\Gamma_F(s_1)(x)-\Gamma_F(s_2)(x)| & \leq & \frac{\chi|D\wh F_M(0)_{|\R^s}|}{(1-\chi\rho)|D\wh F_M(0)_{|\R^u}|}||s_1-s_2||\\
& \leq & \frac{\chi}{1-\chi\rho}\left( \sqrt{\frac{\lambda}{\sigma}}\right)^n||s_1-s_2||.
\end{eqnarray*}
The conclusion of the argument is analogous to the hyperbolic case.
\end{proof}

Using the inverse kergodic maps, we have similar result for the stable vector field. The graph will be the graph of the map $g_M^s$.
As usually, we define 

$$W^s(M)=\disp\bigcup_{n\geq 0} F^{-n}_{( F^n(M))}(W^s_{1}( F^n(M)))
\mbox{ and } W^u(M)=\disp\bigcup_{n\geq 0,}F^{n}(W^u_{1}(F^{-n}(M))).$$
These two manifolds are immersed  Riemannian manifolds. They inherit a Riemannian structure and a Riemannian metric, from the ambient $\R^2$. These two metrics will respectively be denoted by $d^s$ and $d^u$.

For $M$ in $\Lambda\setminus\Lambda_F$ we may define the smallest positive integer  $n$ such that $M_{-n}=f^{-n}(M)$ belongs to $\Lambda_F$. Then we set $W^u_{\rho}(M)$ as the open ball in the Riemannian manifold $f^n(W^u(f^{-n}(M)))$ of radius $\rho$ (and for the metric $d^u$). Similarly we can define $W^s_{\rho}(M)$ for $M$ in $\Lambda\setminus\Lambda_F$. 

\begin{remark} \label{toto}
The result of proposition \ref{prop-vect-loc-inte} and the definitions of the stable or unstable leaves can be adapted to all the points in $\Lambda_{F}$ outside the orbit of tangency. In fact, points in the bottom side of the square have the same stable manifold as $(0,0)$, and points in the left-hand side of the square have the same unstable manifold as $(0,0)$.  
\end{remark}

From now on, we consider that $C_3:=\rho_1.C_0$ is  fixed.
For the next proposition, we give new bounds for $\lambda$ and $\sigma$, respecting the restrictions we put before.

\begin{lemma}\label{invcurves}
There exist $\eps_1 >0$, $0<\lambda<1$ and $\sigma>1$ (satisfying the previous conditions) such that, if $(g(y),y)$ is a ${\CC}^2$ curve contained in region $R_4$, ${\CC}^2$ $\eps_1$-vertical (meaning that the first and second derivatives of $g$ are smaller than $\eps_1$ in modulus), then any curve $(g_0(y),y)$ contained in $f^n(g(y),y)\cap R_4$ is also ${\CC}^2$ $\eps_1$-vertical. 
\end{lemma}

\begin{proof}
First consider a ${\CC}^2$ curve $(s, \phi(s))$, $s\in [-1,1]$, such that $\phi '(0)=0$, $-1<\phi (0)=y_0 \leq 0$ and $2c-1<\phi ''(s)<2c+1$. Then $\phi$ is a convex curve, with an unique critical point $s=0$, admiting two branches of inverse, say $g_1$ and $g_2$ satisfying $g_i(\phi (s))=s$, $i=1,2$. Assume that $g_1(y)=s\geq 0$ and $g_2(y)=s\leq 0$. Then we have

\begin{equation*}
(2c-1)s<\phi '(s)<(2c+1)s
\end{equation*}
and
\begin{equation*}
(2c-1)\frac{s^2}{2}<\phi (s) -\phi (0)<(2c+1)\frac{s^2}{2}
\end{equation*}
Using the last estimate, we find

\begin{equation*}
\sqrt{\frac{\phi(s)-\phi(0)}{2c+1}} <|s|<\sqrt{\frac{\phi(s)-\phi(0)}{2c-1}}
\end{equation*}
We do now estimates for $g_1$, and the ones for $g_2$ are analogous. Remember that $g_1$ is the positive branch of inverse of $\phi$. Then we have

\begin{equation}\label{g1linhamenor}
|g_1'(y)|=\frac{1}{|\phi '(s)|}<\frac{1}{(2c-1)s}<\frac{\sqrt{2c+1}}{(2c-1)\sqrt{y-y_0}}< \frac{\chi_2}{\sqrt{y-y_0}}. 
\end{equation}
We also have
\begin{equation}\label{g1duaslinhas}
|g_1''(y)|=|(g_1'(y))^3\phi ''(s)|<
\frac{\chi_3}{(\sqrt{y-y_0})^3}.
\end{equation}

Now consider the curve $(\lambda^ng_1(y), \sigma^ny)$ for $y$ in the domain of $g_1$ and $1/3<\sigma^ny<1$. This curve is the graph of the function $\psi(y)=\lambda^n g_1(\sigma^{-n}y)$, with $1/3<y<1$, that satisfies (using (\ref{g1linhamenor}))

\begin{equation}\label{psilinha}
|\psi '(y)|=|\lambda^n g_1'(\sigma^{-n}y)\sigma^{-n}|< \frac{\chi_2 \lambda^n}{\sqrt{3}\sigma^{n/2}}
\end{equation}
Using (\ref{g1duaslinhas}), we also find that

\begin{equation}\label{psiduaslinhas}
|\psi ''(y)|=|\lambda^n \sigma^{-2n}g_1''(\sigma^{-n}y)|< \frac{\chi_3 \lambda^n}{\sqrt{3}^3\sigma^{7n/2}}
\end{equation}
Now we conclude the proof of the lemma. Remember that the function $\Gamma$ in remark (\ref{composta}) sends vertical lines into parabolas of equation $y=c(x-q)^2-x_0$. We choose $\eps_1$ such that the image of a ${\CC}^2$ $\eps_1$-vertical curve is sent by $\Gamma$ in a curve $(s,\phi (s))$ satisfying $2c-1<\phi ''(s)<2c+1$. We choose $\lambda$ small  and $\sigma$ big enough such that the last terms in estimates (\ref{psilinha}) and (\ref{psiduaslinhas}) are both smaller than $\eps_1$ for $n=1$.

Start with a ${\CC}^2$ $\eps_1$-vertical curve. Since the linear map $L$ applied to a ${\CC}^2$ $\eps_1$-vertical curve is yet a ${\CC}^2$ $\eps_1$-vertical curve, the map $\Gamma$ guarantees the estimates for the second derivative of the image $(s,\phi (s))$ of the initial curve. Then, up to a change of the coordinate $s$, we can assume that the only critical point of $\phi$ is $0$, and its image is $\leq 0$ (the image of the function $\Gamma$ is contained in the region below $y=c(x-q)^2$). This completes the proof of the lemma. 
\end{proof}

Notice that for a given $c$ and for a given $\eps_1$ the previous
computations give lower bounds for $\sigma$ and $\frac1\lambda$. 
From now on the
constants $c$ and $\eps_1$ are fixed, and $\lambda$ and $\sigma$  are
assumed to satisfy all the required previous conditions. 

\begin{proposition}\label{prop-f-mix}
Let $U$ be an open set in $\Lambda$. Then there exist two positive integers
$n^-_U$ and $n^+_U$ such that $f^{n^+_U}(U)$ contains an unstable
manifold which joins the bottom of the square to the top of the
square, and $f^{-n^-_U}(U)$ contains a stable
manifold which joins the left of the square to the right of the
square. 
\end{proposition}

\begin{proof}
We first check that there exists a point in $U$ whose forward orbit meets
$A$. Indeed, pick any $M$ in $U$ and iterate $W^u_{1}(M)\cap U$. As
long as this connected piece of unstable leaf stays in $R_3$ or $R_5$ its length
growths exponentially fast. Then one of its forward iterates must be a
vertical curve which joints the top of $R_3$ to the bottom of $R_3$ or
the top of $R_5$ to the bottom of $R_5$. Hence the next iterate crosses
$R_4$ and the after next iterate meets $A$. Now, the first iterate of
this last piece of unstable leaf lies in $R'_1$ and joins the bottom
of $\CQ$ to the top of $\CQ$. The same holds for the backward iterates
of $U$, exchanging the vertical and the horizontal directions.

The consequence is that for every $n\geq n^+_U$, $f^n(U)$ contains a
piece of unstable manifold which joins the bottom of $\CQ$ to the top
of $\CQ$. In the same way, for every $m\geq n^-_U$, $f^{-m}(U)$
contains a piece of stable manifold which joins the left side of $\CQ$
to the right side. 
\end{proof}

\begin{cor}
The map $f$ is mixing.
\end{cor}
\begin{proof}
Let $V$ be any open set in $\CQ$. Then for every
$n\geq n^+_U+n^-_V$, $f^{n}(U)\cap V$ is non-empty ; thus $f$ is mixing.
\end{proof}

\begin{remark} \label{rem-struct-prod}
In fact lemma \ref{invcurves}  and proposition \ref{prop-f-mix} prove that the unstable manifolds are long curves, in the sense that for 
any point in $\Lambda^\#$, the connected component of its unstable manifold in $\CQ$ which 
contains the point is a curve which goes from the bottom to the top of
the square. This connected component is denoted by $W^u_{loc}(M)$
Analogously, we define $W^s_{loc}(M)$. Notice that, since the tangent vectors of those curves are inside the cone fields, there exists a positive constant $K$ such that the arc length of $W^{s,u}_{loc}(M)$ is smaller than $K$, for all $M\in \Lambda$.  
As a consequence of the results in the next section, given two manifolds $W^s_{loc}(M)$ and  $W^u_{loc}(M')$, their intersection consists of exactly one point, that is,
$\Lambda^\#$ has a product structure. For two points $M$
and $M'$ in ${\CQ}^\#$ we set $\llb M,M'\rrb$ the intersection of the stable leaf which contains $M$ in the square and the unstable leaf which contains $M'$.
\end{remark}

\subsection{Extra properties of the foliations and $F$}

In the next result we use the existence of (un)stable manifolds to  prove the H\"older regularity for the hyperbolic splitting.
\begin{proposition}\label{lem-vect-hold}
There exist two positive constants, $C_1$\index{$C_1$} and
$C_2$\index{$C_2$} such that in $\Lambda_F$, the maps $z\mapsto e^u(z)$ and
$z\mapsto e^s(z)$ are respectively $\disp {C_1}-\frac12$-H\"older
continuous and $\disp C_2-\frac12$-H\"older continuous for the Euclidean metric $||.||$.
\end{proposition}

\begin{proof}
Let us first consider two $C^2$ maps,
$\phi_1:[0,\frac13]\rightarrow [0,1]$ and $\phi_2:[0,\frac13]\rightarrow
[0,1]$ with first and second
derivatives smaller, in modulus, than $\epsilon>0$. Assume that the two graphs have  empty intersection, and take $I=(a,A)$ in the graph of $\phi_1$ and $J=(b,\phi_2(b))$ in the graph of
$\phi_2$. Let $K=(a,\phi_2(a))$, $B:=\phi_2(a)$,
and assume that $B\geq A$, $\phi_1'(a)\geq
\phi_2'(a)$ and $a\leq \frac16$ (the other cases can be obtained from this one). We thus have 
\begin{equation}\label{equ1-hold-eu}
|\phi_1'(a)-\phi_2'(b)|\leq
 |\phi_1'(a)-\phi_2'(a)|+|\phi_2'(a)-\phi_2'(b)|.
\end{equation}
The conditions on $\phi_2$  mean that the second term in the right hand
side of (\ref{equ1-hold-eu}) is smaller than $\epsilon.|a-b|$. Now, the standard
computation in differential calculus gives 
$$\phi_1(t)\geq A+(t-a)\phi_1'(a)-\epsilon.\frac{(t-a)^2}{2},\mbox{ and }
\phi_2(t)\leq B+(t-a)\phi_2'(a)+\epsilon.\frac{(t-a)^2}{2}.$$ 
The fact that the two curves have an empty intersection yields to

$$\epsilon.(t-a)^2+(t-a)(\phi_2'(a)-\phi_1'(a))+B-A\geq
0,$$
for every $t$ in $[0,\frac13]$. The minimal value is obtained for
$\disp t-a:=\frac{\phi_1'(a)-\phi_2'(a)}{2\epsilon}$. \\
If this value is
smaller than $\frac13-a$, then we obtain 

$$\epsilon.\left(\frac{\phi_1'(a)-\phi_2'(a)}{2\epsilon}\right)^2-\frac{(\phi_1'(a)-\phi_2'(a))^2}{2\epsilon}+B-A\geq
0,$$
which yields to

\begin{equation}\label{equ2-hold-eu}
0\leq (\phi_1'(a)-\phi_2'(a))^2\leq 4\epsilon.(B-A).
\end{equation}
If $$\disp \frac{\phi_1'(a)-\phi_2'(a)}{2\epsilon}\geq \frac13-a,$$ then we
obtain $$\disp
\epsilon.\left(\frac13-a\right)^2-\left(\frac13-a\right)(\phi_1'(a)-\phi_2'(a))+B-A\geq
0,$$ which yields to

\begin{equation}\label{equ3-hold-eu}
0\leq \phi_1'(a)-\phi_2'(a)\leq 12(B-A).
\end{equation}
Hence, (\ref{equ1-hold-eu}) (\ref{equ2-hold-eu}) and
(\ref{equ3-hold-eu}) yield to

$$|\phi_1'(a)-\phi_2'(b)|\leq C|a-b|+C'\sqrt{B-A}.$$
We also have $$\disp \sqrt{B-A}\leq
\sqrt{|\phi_2(b)-\phi_1(a)|}+\sqrt{|\phi_2(b)-\phi_2(a)|}.$$ This last
term in the sum is smaller than $\disp\sqrt\epsilon.|a-b|^\frac12$, which
is smaller than $\disp\sqrt{\epsilon}||I-J||^\frac12$. 

As a direct consequence of our computation, with $\epsilon=\eps_1$, we get that the unstable
vector field is at least $\frac12$-H\"older continuous in $R_3\cup
R_4\cup R_5$, with some universal constant (which depends on $\eps_1$). Now the
map $\Gamma$ is $C^2$ on the compact set $[0,1]^2$, and so there
exists some constant $C_1$ such that the unstable vector field is (at
least) $\frac12$-H\"older continuous in $\Lambda_F$. The proof
for the stable vector field is analogous. 
\end{proof}


\section{Markov partitions and equilibrium states}\label{sec-thermo}
In this section we prove that the map $f$ is not expansive (see the definition below), and construct a H/"older continuous finite-to-one semi-conjugacy from the full shift on the space of sequences of three symbols to $\Lambda$.

\subsection{Expansiveness fails}

Recall that a diffeomorphism $g$ is said to be expansive if there exists a constant $\delta>0$ such that, if $d(g^n(x),g^n(y))<\delta$ for every $n\in \Z$, then $x=y$.
Expansiveness is the usual tool to get existence of the equilibrium states. Indeed, expansiveness implies the upper semi-continuity for the map $\mu\mapsto h_\mu(g)$, where $h_\mu(g)$ is the metric entropy of $g$ with respect to $\mu$. Nevertheless, expansiveness is not necessary to get this upper semi-continuity. In Bowen's proof (proposition 2.19 p. 64 in
\cite{bowen}), the main argument to get it is  that every  partition with sufficiently small radius is generating (its entropy is equal to  the metric entropy). In our case, it is not true, due to the non-expansiveness, that every partition with  sufficiently small diameter is generating. However it can be proved that the geometric partition $\CG_1$, to be defined below, is generating. This is sufficient to get the upper semi-continuity of the metric entropy.

\begin{figure}[htb]
\begin{center}
\includegraphics[width=2in]{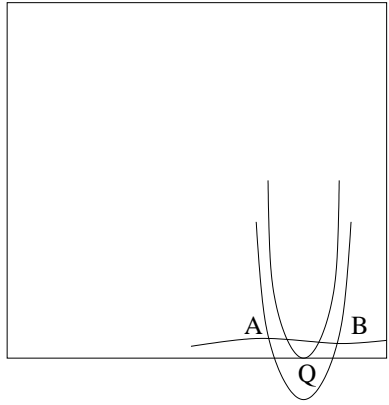}	
	\caption{	\label{nonexp} choice of $A$ and $B$}
\end{center}
\end{figure}

Let us now briefly prove why $f$ is not expansive. Let  $\delta>0$ be fixed, and pick two points $A$ and $B$  close  to the critical point $Q=(q,0)$, such that $B$ belongs to the unstable manifolds of $A$ but just on the opposite side of $Q$, and $B$ also belongs to the stable manifold of $A$, see figure \ref{nonexp}. This can be done because, locally, the unstable manifolds are $\CC^2$-close of two branches of parabolas and the stable manifolds are long ``horizontal'' curves. Choosing $A$ and $B$ sufficiently close to $Q$ we  get that $A\neq B$ but $d(f^n(A),f^n(B))\leq \delta$ for every $n\in\Z$.

\subsection{Proof of Theorem B} 
  
  For $n\in \N$, we define the degenerated partition $\CG_n$ as follows.
  Let $\CR^i_n$, $i\in\{ 1,\ldots,3^{2n}\}$ be the connected components of 
  $f^n(\CQ^\#)\cap f^{-n}(\CQ^\#)$. Then $\CG_n= \{ \overline{\CR^i_n}, i=1,\ldots,3^{2n}\}$.
  For convenience, we call $\CR^1_n$ the element of $\CG_n$ to which $(0,0)$ belongs. Then $\CG_n$ 
  has the following properties:
  
  \begin{itemize}
  \item[-]each atom $\CR^i_n$ is bounded by arcs of the stable and unstable manifolds of 
  $(0,0)$ and $(1,1)$;
  \item[-]for each $n\in\N$ and each $i\in\{ 1,\ldots,3^{2n}\}$, $f^n(\CR^i_n)$ is a stripe crossing $\CQ$ 
  from the bottom to the top, and $f^{-n}(\CR^i_n)$ is a stripe crossing $\CQ$ 
  from the left to the right sides;
  \item[-]there are $2n$ pairs of atoms that are not disjoint: they intersect in a point of 
   the orbit of tangency;
  
  \item[-]each atom $\CR^i_{n+1}$ of $\CG_n$ can be obtained as $f^{-1}[f(\CR^j_n)\cap \CR^l_n]\cap
   f[f^{-1}(\CR^j_n)\cap \CR^k_n]$, where $\CR^i_{n+1}\subset \CR^j_n$.  
  \end{itemize}

\begin{figure}[htb]
\begin{center}
\includegraphics[width=4in]{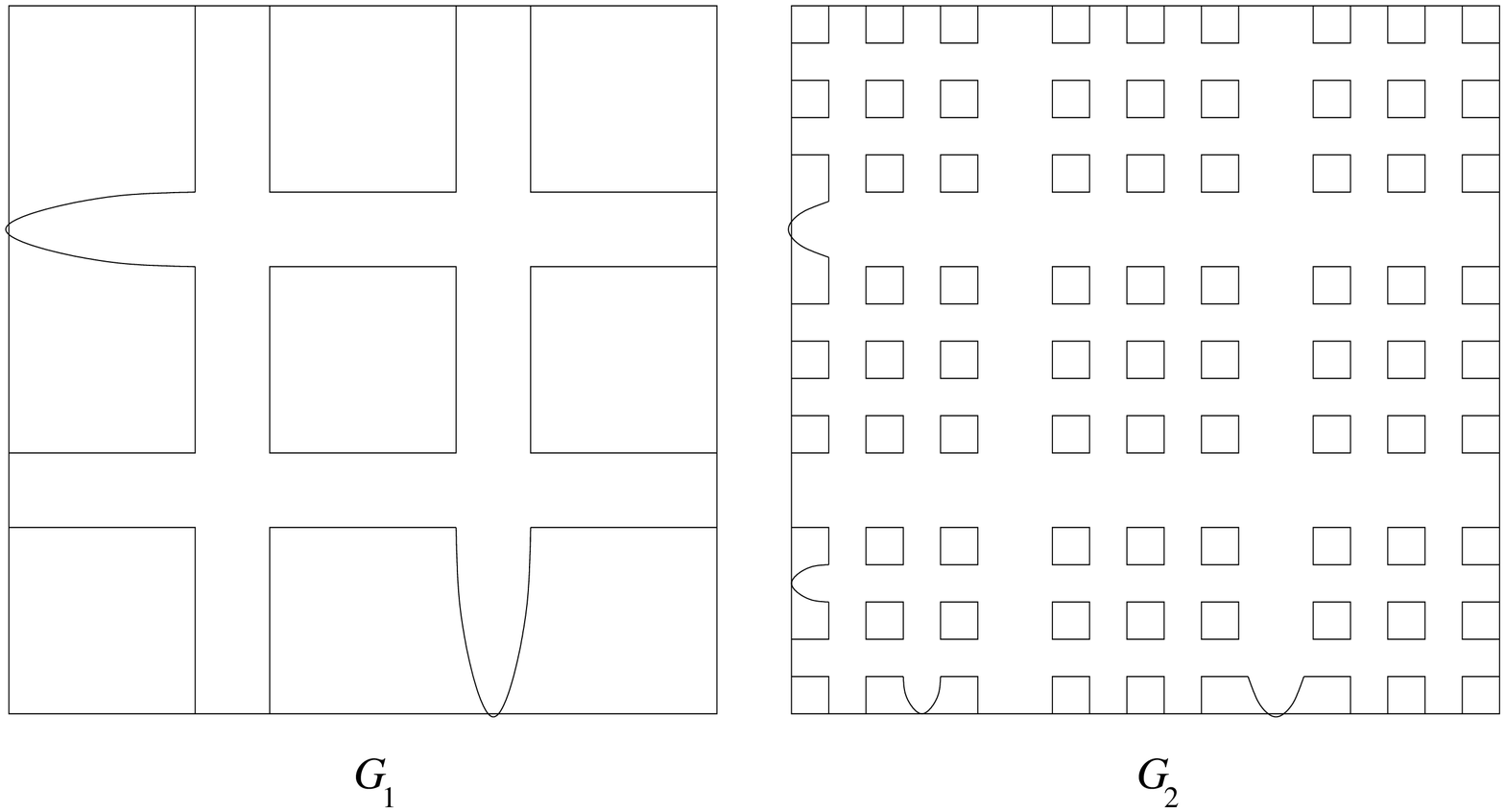}	
	\caption{	\label{partition}the partitions $\CG_1$ and $\CG_2$}
\end{center}
\end{figure}

 \begin{remark}\label{itinerary}
  The partitions $\CG_n$ are not real partitions, in the sense that some points of the orbit of
  homoclinic tangency are in two atoms. Due to the last condition above, for each $-n+1\leq k \leq n-1$,
  the image by $f^k$ of points 
  in the same element of $\CG_n$ belong to the same element of $\CG_1$. 
\end{remark}

\begin{proposition}\label{atoms}
 The length of the boundary of each $\CR^i_n$ goes to zero exponentially with $n$.  
\end{proposition}
\begin{proof}
Let $\CR = \CR_n^i\in \CG_n$. Denote by $\partial ^u\CR$ and $\partial ^s\CR$, respectively, the 
unstable and stable components of the boundary of $\CR$, each containing two connected components.

Notice that, due to remark \ref{itinerary}, we can define the itinerary of $\CR$ from $-n+1$ to $n-1$. 
For some values of $k$ in such interval, $f^k(\CR)$ may be contained in $\CR^1_1$. Define $n_i$ and 
$n_{-i}$ as follows: $n_1=k-1$ where $k$ is the first iterate of $\CR$ in $\CR^1_1$, 
and $n_{i+1}>n_{i}$ is the biggest iterate of $\CR$ not in $\CR^1_1$ and such that 
$f^{n_{i+1}}(\CR)\subset \CR^1_1$. The sequence $n_{-i}$ is defined analogously for iterates by $f^{-1}$.

There are some cases to be considered concerning the itinerary of $\CR$. The easiest 
one is when the iterates $-n+1$ and $n+1$ are outside $\CR^1_1$, and are not at the sequences $n_i$ and $n_{-i}$.
 In this case, all the unstable vectors $v_u$ tangent to $\partial^u\CR$, and the stable vectors $v_s$, tangent to $\partial^s\CR$, have grown exponentially by $Df^n$ and $Df^{-n}$, with a factor of $\sqrt{\sigma}$ and $\sqrt{\lambda}$, respectively. 
 Then we have $||Df^n_Mv||\geq \sigma^{n/2}$ and $||Df^{-n}_{M'}v_s||\geq \lambda^{-n/2}$, for $M$ in $\partial^u\CR$ and $M'$ in $\partial^s\CR$. Since the length of the local stable and unstable manifolds in $\CQ$ is bounded from above by  $K$ (see remark \ref{rem-struct-prod}), we have that each component of $\partial^u\CR$ has length smaller than $K\sigma^{-n/2}$, and each component of $\partial^s\CR$ has length smaller than $K\lambda^{n/2}$.
 
 We now prove the proposition in the other cases for the unstable boundary of $\CR$, and the stable case is completely analogous. Assume first that $n_i=n-1$. It means that $f^{n-1}\CR$ is contained in one of the two atoms of $\CG_1$ that contain $(q,0)$, and its length at that moment is smaller than $K'$. Since the estimates above are valid until $n-1$, we have that $\partial^u\CR$ has length smaller than $K'\sigma^{-(n-1)/2}<\tilde{K}\sigma^{-n/2}$.
 
 In the remaining case, we have that for each $k$ between $n_i$ and $n-1$, $f^k(\CR)\subset \CR^1_1$. Notice that, in this case, we have that the estimates on the growth of the tangent vector to  $\partial^u\CR$ apply until $n_i$. There is also a part of  $f^{n_i}\partial^u\CR$ that escapes from $R^1$ by the remaining iterates, assuming size smaller than $K$, to which our estimates apply until $n$. The remaining part of  $f^{n_i}(\partial^u\CR)$ is confined to the horizontal stripe $[0,1]\times[0,\sigma^{-(n-n_i)}]$. The maximum length of an unstable manifold in that stripe is smaller than $\chi\sqrt{\sigma^{n-n_i}/c}$. It gives us that the length of  $\partial^u\CR$ is smaller than $K\sigma^{-n/2}+\chi'\sqrt{\sigma^{n-n_i}/c}\sigma^{-n_i/2}=\tilde{\tilde{K}}\sigma^{-n/2}$.
\end{proof}

Let $\S_3$ be the space of bi-infinite sequences of three symbols. From now on, we also denote by $\sigma$ the shift on $\S_3$, avoiding confusion with the unstable eigenvalue, would it appear. As a corollary of proposition \ref{atoms}, we get the next important result:

\begin{proposition}
There exists a finite-to-one and onto  H\"older continuous semi-conjugacy, $\Theta$, between the two dynamical systems $(\S_3,\sigma)$ and $(\Lambda,f)$.
\end{proposition}
\begin{proof}
We associate  the numbers 0,1, and 2 to the components of $f(\CQ)\cap\CQ$ as follows. We assign the number $0$  to $R_1'$, $1$ to the component at the right-hand side of the critical point $(q,0)$, and  $2$ to the remaining part. Notice that the point $(q,0)$ is in the intersection of the two regions $1$ and $2$. Via this correspondence, we associate to each atom $\CR_n^i$ of the partition $\CG_n$, a {\it centered word} of length $2n+1$, in the following way. As a consequence of the definition of $\CG_n$, all points in  $\CR_n^i$ have their images by $f^k$ falling in the same stripe $s_k=0$, $1$ or $2$, for all $-n\leq k\leq n$. We associate to $\CR_n^i$ the centered word $[s_{-n},\ldots,s_n]$. Notice that 2 different atoms of $\CG_n$ define two different words.

Let us now define the map $\Theta:\S_3\rightarrow \Lambda$. Let  $\ul{\xi}\pardef(\xi_k)_k$ be in $\S_3$. As usually we define the $n$-cylinder $C_n(\ul\xi)$ as the set of sequences  $\ul{\xi'}=(\xi_k')_k$ such that $\xi_k=\xi_k'$ for any $-n\leq k\leq n$. Such a $n$-cylinder defines (and is defined by) a  centered word of length $2n+1$ in $\S_3$.
To each such $n$-cylinder we associate the unique element of the partition $\CG_n$ which have the same centered word: we  set $\Theta(\CC)=R$, where $\CC$ is a $n$-cylinder, and $R$ is the atom of $\CG_n$ which have the same centered word  of length $2n+1$ than $\CC$. For a fixed $\ul\xi$ in $\S_3$, each $\Theta(C_n(\ul\xi))$ is a non-empty compact set ; the sequence  $(\Theta(C_n(\ul\xi)))_n$
 is a decreasing sequence of compact sets, which thus have a non-empty intersection. Proposition \ref{atoms} also implies that this infinite intersection is a single point. Hence we set 
$$\Theta(\ul\xi)\pardef\bigcap_n\Theta(C_n(\ul\xi)).$$
This defines the map $\Theta$. Clearly we have $\Theta\circ \s=f\circ \Theta$.

We now prove that $\Theta$ is onto. 
For any $\xi$ in $\Lambda$ we build the code $\ul\xi$ in the following way. For each integer $k$, $\xi_k$ is the number of the (one) full vertical band which contains $f^k(\xi)$. This defines at least one bi-infinite $\ul{\xi}\pardef(\xi_k)_k$ in $\S_3$ such that $\Theta(\ul\xi)=\xi$.

To see that $\Theta$ is  finite-to-one, notice that, by construction of the map, only points in the critical orbit have several pre-images by $\Theta$. Indeed  the atoms of the partition $\CG_n$ are disjoint, except for the points in the critical orbit.
Moreover, any such point $\xi$ as at least two pre-images by $\Theta$ which are the bi-infinite words $\ldots,0,0,1,0,0,\ldots$ and $\ldots,0,0,2,0,0,\ldots$ where the 1 and the 2 are at the same position in the sequence. However itinerary of the critical orbit is very simple. Any point $\xi$ whose orbit contains $T$ (and $Q$), $f^k(\xi)$ has its backward iterates the segment $\{0\}\times [0,t]$, its first iterate is $Q$ and the remaining iterates are in  $[0,q]\times\{0\}$. As the `0' band is disjoint from the 2 other bands, a critical point cannot have other pre-images than the 2 ones  indicated above. This completes the prove that $\Theta$ is finite-to-one.

To conclude the proof of the proposition, it remains to show that $\Theta$ is H\"older continuous. It follows from  classic arguments, that is based in the exponential decay for the diameter of the atoms of $\CG_n$.  We recall that the distance on $\S_3$ between two sequences $\ul\xi$ and $\ul\xi'$ is $\disp \frac{1}{2^n}$, where $n$ is the largest non-negative integer such that $\ul\xi'$ belongs to $C_n(\ul\xi)$. As a consequence of the proof of proposition \ref{atoms}, for any non-negative integer $n$ and for any $\ul\xi$ in $\S_3$,  $\Theta(C_n(\ul\xi))$ has  diameter smaller than $K.(\lambda^{\frac{n}{2}}+\sigma^{-\frac{n}{2}})$. If $\ul\xi'\in C_n(\ul\xi)$ but $\ul\xi'\notin C_{n+1}(\ul\xi)$, then $d(\ul\xi,\ul\xi')=1/2^n$, and $\Theta(\ul\xi)$ and $\Theta(\ul\xi')$ belong to the same atom of $\CG_n$.   Therefore we get

$$\left|\Theta(\ul\xi)-\Theta(\ul\xi')\right|\leq 2.K. d^\gamma(\ul\xi,\ul\xi'),$$
where $\disp\gamma=\min(-\frac{\log\sqrt{\lambda}}{\log2},\frac{\log\sqrt{\sigma}}{\log2})$.

\end{proof}

\begin{remark}\label{rem-orbi-crit}
Notice that for every $f$-invariant measure $\mu$, the critical orbit has null $\mu$-measure. In the same way, for every $\s$-invariant measure in $\S_3$, the set of sequences where $\Theta$ is one-to-one has full measure. Therefore, any push-forward of any $\s$-invariant measure in $\S_3$ on $\Lambda$ by $\Theta$ is a $f$-invariant measure; conversely any pull-back on $\S_3$ of any  $f$-invariant measure is a $\s$-invariant measure. Notice also that, as mentioned in remark \ref{rem-struct-prod}, if $M$ and $M'$ are not in the orbit of tangency, $W^s_{loc}(M)\cap W^u_{loc}(M')$ consists of only one point, that we call $\llb M,M'\rrb$. We also have that each $\CR^i_n$ is a rectangle, in the sense of Bowen: if $M$ and $M'$ are in $\CR^i_n$ then $\llb M,M'\rrb \in \CR^i_n$; 
\end{remark}

The proof of Theorem B is now immediate. Indeed, for any H\"older continuous function $\varphi:\CQ\rightarrow \R$, $\varphi\circ\Theta$ is a H\"older continuous function from $\S_3$ to $\R$. Thus there exists a unique equilibrium state, which is also a Gibbs measure, associated to $\varphi\circ\Theta$ in $\S_3$. The push-forward of this measure on $\Lambda$ gives some $f$-invariant measure $\mu_\varphi$ on $\Lambda$. As $\Theta$ is finite-to-one, $\mu_\varphi$ has maximal $\varphi$-pressure on $\Lambda$.


\bibliographystyle{plain}
\bibliography{mabiblio}
\end{document}